\documentclass{amsart}

\setbox0=\hbox{$+$}
\newdimen\plusheight
\plusheight=\ht0
\def\+{\;\lower\plusheight\hbox{$+$}\;}

\setbox0=\hbox{$-$}
\newdimen\minusheight
\minusheight=\ht0
\def\-{\;\lower\minusheight\hbox{$-$}\;}

\setbox0=\hbox{$\cdots$}
\newdimen\cdotsheight
\cdotsheight=\plusheight
\def\cds{\lower\cdotsheight\hbox{$\cdots$}}

\newcommand{\df}{\dfrac}

 \renewcommand{\a}{\alpha}

\renewcommand{\t}{\varphi}

\renewcommand{\(}{\left\(}
\renewcommand{\)}{\right\)}

\renewcommand{\i}{\infty}
\renewcommand{\b}{\beta}
\renewcommand{\pmod}[1]{\,(\textup{mod}\,#1)}

\def\receivedline{\relax}
\def\dedication#1{\receivedline\vskip4pt
{\normalsize\begin{center}#1\end{center}\vskip1sp}}

\newcommand{\beqs}{\begin{equation*}}
\newcommand{\eeqs}{\end{equation*}}
\numberwithin{equation}{section}
 \theoremstyle{plain}
\newtheorem{theorem}{Theorem}[section]

\newtheorem{corollary}[theorem]{Corollary}

\newtheorem{remark}[theorem]{Remark}

\begin{document}

\title[Quadratic forms ] {Ramanujan's identities and
representation of integers by certain binary and quaternary
quadratic forms}
\author{Alexander Berkovich}
\author{Hamza Yesilyurt}

\address{Department of Mathematics, University of Florida, 358 Little Hall,  Gainesville, FL
  32611, USA}\email{alexb@math.ufl.edu}

\email{hamza@math.ufl.edu} \keywords{quadratic forms, $q$-series
identities, $eta$-quotients, multiplicative functions}
\subjclass[2000]{Primary: 11E16, 11E25, 11F27, 11F30; Secondary:
05A19, 05A30, 11R29}
\thanks{This research is supported in part by NSA Grant MSPF-06G-150.}

\begin{abstract}

We revisit old conjectures of Fermat and Euler regarding
representation of integers by binary quadratic form  $x^2+5y^2$.
Making use of  Ramanujan's  $_1\psi_1$ summation formula we
establish a new Lambert series identity for $\sum_{n,m=-\infty}^{\infty}
q^{n^2+5m^2}$. Conjectures of Fermat and Euler  are shown to
follow easily from this new formula. But we don't stop there.
Employing various formulas found in Ramanujan's notebooks and
using a bit of ingenuity we obtain a collection of new Lambert
series for certain infinite products associated with quadratic
forms such as $x^2+6y^2$, $2x^2+3y^2$,  $x^2+15y^2$, $3x^2+5y^2$,
$x^2+27y^2$, $x^2+5(y^2+ z^2+ w^2)$, $5x^2+y^2+ z^2+ w^2$.
In the process, we find many new multiplicative $eta$-quotients and
determine their coefficients.

\end{abstract}

\maketitle
 \dedication{  }
\section{Introduction}

A binary quadratic form (BQF) is a function
\begin{equation*}
Q(x,y)=ax^2+bxy+cy^2
\end{equation*}
with $a, b, c$ $\in$ $\mathbb{Z}$. It will be denoted by $(a,b,c)$.
We say that $n$ is represented by $(a,b,c)$ if there exist $x$ and
$y$ $\in$ $\mathbb{Z}$ such that $Q(x,y)=n$.

The representation theory of BQF has a long history that goes back
to antiquity. Diophantus' Arithmeticae contains the following
important example of composition of two forms

\begin{equation*}
(x_1^2+y_1^2)(x_2^2+y_2^2)=(x_1x_2-y_1y_2)^2+(x_1y_2+x_2y_1)^2.
\end{equation*}

Influenced by Diophantus, Fermat studied representations by
$(1,0,a)$. For $a=1,2,3$ he proved a number of important results
such as the following.
 \bigskip

\textit{A prime $p$ can be written as a sum of two squares iff
$p\equiv 1 \pmod{4}$}.
 \bigskip

We remark that representation by $(1,0,3)$ played an important
role in Euler's proof of Fermat's Last Theorem in the case of
$n=3$.

Fermat realized that $(1,0,5)$ was very different from the
previous cases $(1,0,1)$, $(1,0,2)$ and $(1,0,3)$ considered by
him. He made the following conjecture.
 \bigskip

\textit{If $p$ and $q$ are two primes that are congruent to $3$ or $7$
$\pmod{20}$, then $pq$ is representable by $(1,0,5)$}.
 \bigskip

Euler made two conjectures that were very similar to those of
Fermat:

 \bigskip

\textit{a.  Prime $p$ is represented by $(1,0,5)$ iff   $p \equiv
1\,\text{or}\; 9 \pmod {20}$}.

 \bigskip
 \textit{b. If $p$ is prime then $2p$ is represented by
$(1,0,5)$ iff $p \equiv 3 \,\text{or}\, 7 \pmod {20}$}.
 \bigskip

However, his next conjecture for $(1,0,27)$ contained an unexpected
cubic residue  condition:

 \bigskip

\textit{ Prime $p$ is represented by $(1,0,27)$ iff  $p \equiv 1
\pmod 3$ and $2$ is a cubic residue modulo $p$}.

 \bigskip

Lagrange and Legendre initiated systematic study of quadratic
forms. But it was Gauss who brought the theory of BQF to
essentially its modern state. He introduced class form groups and
genus theory for BQF. He proved Euler's conjecture for $(1,0,27)$
and in the process discovered a so-called cubic reciprocity law.
Gauss' work makes it clear why $(1,0,27)$  is much harder to deal
with than $(1,0,5)$. Indeed, a class form group with discriminant
$-20$ consist of two inequivalent classes $(1,0,5)$ and $(2,2,3)$.
These forms can't represent the same integer. On the other hand, a
class form group with discriminant $-108$ consist of three classes
$(1,0,27)$, $(4,2,7)$, $(4,-2,7)$. These forms belong to the same
genus. That is, they may represent the same integer. An interested
reader may want to consult \cite{DC} and \cite{mw} for the wealth
of historical information and \cite{ws} for the latest
development.

In his recent book, Number Theory in the Spirit of Ramanujan,
Bruce Berndt discusses representation problem for $(1,0,1)$,
$(1,0,2)$, $(1,1,1)$, $(1,0,3)$. Central to this approach is
Ramanujan's $_1\psi_1$ summation formula which implies in
particular that \cite[p.58, eq.~(3.2 .90)]{berndt2}
\begin{equation}\label{iiin}
\sum_{x,y \in \mathbb{Z}}q^{x^2+y^2}=1+4\sum_{n \geq
1}\df{q^n}{1+q^{2n}}.
\end{equation}

Using geometric series it is straightforward to write the right
hand side of \eqref{iiin} as
\begin{align*}
1+4\sum_{n \geq 1, m\geq 0}(-1)^mq^nq^{2nm}&=1+4\sum_{n \geq 1,
m\geq 0}(-1)^mq^nq^{n(2m+1)}\\&=1+4\sum_{n \geq 1, m\geq
1}\Bigr(\df{-4}{m}\Bigl)q^{nm}\\&=1+4\sum_{n \geq 1
}\sum_{d|n}\Bigr(\df{-4}{d}\Bigl)q^{n},
\end{align*}

where we used the Kronecker symbol to be defined in the next
section and the well known formula
\begin{equation*}
\Bigr(\df{-4}{n}\Bigl)=\left\{ \begin{array}{ll}
         0  &  \text{if}\;\; n \;\text{is even}, \\
         (-1)^{(n-1)/2}  &  \text{if}\;\; n\; \text{is odd}. \end{array} \right.
         \end{equation*}

Let $r(n)$ be the number of representations of a positive integer
$n$ by  quadratic form $k^2+l^2$. Suppose the prime factorization
of $n$ is given by
\begin{equation*}
n=2^a\prod_{i=1}^{r}p_i^{v_i}\prod_{j=1}^sq_j^{w_j},
\end{equation*}
where  $p_i \equiv 1 \pmod {4}$ and $q_i \equiv -1 \pmod {4}$.

Using the fact that $\sum_{d|n}\Bigr(\df{-4}{d}\Bigl)$ is
multiplicative, we find that

\begin{equation}
r(n)=4 \prod_{i=1}^{r} (1+v_i) \prod_{j=1}^s\df{1+
(-1)^{w_j}}{2}.\label{sofs}
\end{equation}

Reader may wish to consult \cite{apostol} for background on
multiplicative functions, convolution of multiplicative functions
and Legendre's symbol. Clearly, Fermat's Theorem is an immediate
corollary of \eqref{sofs}.

The main object of this manuscript is to reveal new and exciting
connections between
 the work of Ramanujan and the theory of  quadratic forms. This paper is
organized as follows.

We collect necessary definitions and formulas in Section 2.

In Section 3 we use the $_1\psi_1$ summation formula to prove new
generalized Lambert series identities for
\begin{equation*}
\sum_{n,m=-\i}^{\i}q^{n^2+5m^2}  \;\;\; \text{and}\;\;\;
\sum_{n,m=-\i}^{\i} q^{2n^2+2nm+ 3m^2}.
\end{equation*}
These results enable us to derive  simple formulas for the number
of representations of an integer $n$ by  $(1,0,5)$ and $(2,2,3)$.
Conjectures of Fermat and Euler for $(1,0,5)$ are easy corollaries
of these formulas. Our treatment of $(1,0,6)$ and $(2,0,3)$ in
Section 4 is very similar. However, in addition to the $_1\psi_1$
summation formula we  need  to use two cubic identities of
Ramanujan. In Section 5 we treat $(1,0,15)$ and $(3,0,5)$. The
surprise here is that we need to employ one of the fourty
identities of Ramanujan for the Rogers-Ramanujan functions.
Section 6 deals with $(1,0,27)$ and $(4,2,7)$.  We do not confine
our discussion solely to BQF. In Section 7 we boldly treat
quaternary forms $x^2+ 5 (y^2+z^2+w^2)$ and $5x^2+ y^2+z^2+w^2$.
We conclude with a brief description of the prospects for future
work.

\bigskip

\section{Definitions and Useful Formulas}
Throughout the manuscript we assume that $q$ is a complex number
with $|q|<1$. We adopt the standard notation
\begin{align*}
(a;q)_n&:=(1-a)(1-aq)\ldots(1-aq^{n-1}),\notag\\
(a;q)_\i&:=\prod_{n=0}^\i(1-aq^n), \notag\\
E(q)&:=(q;q)_\i.
\end{align*}

Next, we  recall  Ramanujan's definition for a general theta
function. Let

\begin{equation}\label{generaltheta}
f(a,b) := \sum_{n=-\i}^{\i}a^{n(n+1)/2}b^{n(n-1)/2}, \qquad |ab| <
1.
\end{equation}
The function $ f(a,b) $ satisfies the well-known Jacobi triple
product identity \cite[p.~35, Entry 19]{III}
\begin{equation}\label{19III}
f(a,b) = (-a;ab)_{\i}(-b;ab)_{\i}(ab;ab)_{\i}.
\end{equation}
Two important special cases of \eqref{generaltheta} are
\begin{equation}\label{22i}
\varphi(q) := f(q,q) = \sum_{n=-\i}^{\i}q^{n^2} =
(-q;q^2)_{\i}^2(q^2;q^2)_{\i}=\df{E^5(q^2)}{E^2(q^4)E^2(q)},
\end{equation}
and
\begin{equation}\label{22ii}
\psi(q) := f(q,q^3) = \sum_{n=-\i}^{\i}q^{2n^2-n}
=(-q;q^4)_\i(-q^3;q^4)_\i(q^4;q^4)_\i=\df{E^2(q^2)}{E(q)} .
\end{equation}
The product representations in \eqref{22i}--\eqref{22ii} are special
cases of \eqref{19III}.
We shall use the famous quintuple product identity, which, in
Ramanujan's notation,  takes the form \cite[p.~80, Entry
28(iv)]{III}
\begin{equation}\label{qtp}
E(q)\df{f(-a^2,-a^{-2}q)}{f(-a,-a^{-1}q)}
=f(-a^3q,-a^{-3}q^2)+af(-a^{-3}q,-a^{3}q^2),
\end{equation}
where $ a $ is any complex number.

Function $ f(a,b) $ also satisfies a useful addition formula.
For each nonnegative integer $ n$, let
\begin{equation*}
U_n := a^{n(n+1)/2}b^{n(n-1)/2} \qquad \text{and} \qquad V_n :=
a^{n(n-1)/2}b^{n(n+1)/2}.
\end{equation*}
Then \cite[p.~48, Entry 31]{III}
\begin{equation}\label{31III}
f(U_1,V_1) =
\sum_{r=0}^{n-1}U_rf\left(\df{U_{n+r}}{U_r},\df{V_{n-r}}{U_r}\right).
\end{equation}

From \eqref{31III} with $n=2$, we obtain
\begin{equation}\label{2dd}
f(a,b)=f(a^3b,ab^3)+af(\df{b}{a},\df{a}{b}(ab)^4).
\end{equation}
A special case of \eqref{2dd} which we frequently use is
\begin{equation}
\t(q)=\t(q^4)+2q\psi(q^8).\label{tdd}
\end{equation}
With $a=b=q$ and $n=3$, we also find that
\begin{equation}
\t(q)=\t(q^9)+2qf(q^3,q^{15}).\label{tdd3}
\end{equation}

Our proofs employ a well-known special case of Ramanujan's $_1\psi_1$ summation formula:\\
If $|q|<|a|<1$, then \cite[p.~32, Entry 17]{III}
\begin{equation}\label{1psi1}
E^3(q)\df{f(-ab,-q/ab)}{f(-a,-q/a)f(-b,-q/b)}=\sum_{n=-\i}^{\i}\df{a^n}{1-bq^n}.
\end{equation}

We frequently use  the elementary result \cite[p.~45, Entry
29]{III}. If $ab=cd$, then
\begin{align}\label{ff}
f(a,b)f(c,d)=f(ac,bd)f(ad,bc)+af(\df{b}{c},\df{c}{b}abcd)f(\df{b}{d},\df{d}{b}abcd).
\end{align}

Next, we recall that for an odd prime $p$, Legendre's Symbol
$\Bigr(\df{n}{p}\Bigl)$ or $(n \mid p)$ is defined by

\begin{equation*}
\Bigr(\df{n}{p}\Bigl)=\left\{ \begin{array}{ll}
         1  &  \text{if n is a quadratic residue modulo}\; p, \\
         -1  & \text{if n is a quadratic nonresidue modulo}\; p . \end{array} \right.
\end{equation*}

Kronecker's Symbol $\Bigr(\df{n}{m}\Bigl)$ is defined as follows\\

\begin{equation*}
\Bigr(\df{n}{k}\Bigl)=\left\{ \begin{array}{ll}
         1  &  \text{if}\;\; k=1, \\
         0  &  \text{if}\;\; k\; \text{is a prime dividing $n$,}\\
         \text{Legendre's symbol}  & \text{if $k$ is an odd prime}. \end{array} \right.
\end{equation*}

\begin{equation*}
\Bigr(\df{n}{2}\Bigl)=\left\{ \begin{array}{ll}
         0  &  \text{if}\;\; n \;\text{is even}, \\
         1  &  \text{if}\;\; n\; \text{is odd},\; n \equiv \pm 1 \pmod8,\\
         -1 & \text{if}\;\; n \;\text{is odd},\; n \equiv \pm 3 \pmod8. \end{array} \right.
\end{equation*}

$$ \text{In general}, \;\;\Bigr(\df{n}{m}\Bigl)=\prod_{i=1}^{s}\Bigr(\df{n}{p_i}\Bigl)
\;\;\text{if}\;\; m=\prod_{i=1}^{s}p_i\;\;\text{is a prime factorization of}\;\; m. $$

It is easy to show that
$\Bigr(\df{a}{bc}\Bigl)=\Bigr(\df{a}{b}\Bigl)\Bigr(\df{a}{c}\Bigl)$
and
$\Bigr(\df{ab}{c}\Bigl)=\Bigr(\df{a}{c}\Bigl)\Bigr(\df{b}{c}\Bigl)$.
Hence, $\Bigr(\df{n}{m}\Bigl)$ is a completely multiplicative
function of $n$ and also of $m$.

\bigskip

\section{Lambert Series Identities for $\sum_{n,m=-\i}^{\i}q^{n^2+5m^2}$}
\begin{theorem}\label{d20}
\begin{align}
\t(q)\t(q^5)&=2\Bigr\{\sum_{n=-\i}^{\i}\df{q^n}{1+q^{10n}}-\sum_{n=-\i}^{\i}\df{q^{5n+2}}{1+q^{10n+4}}\Bigl\}\label{d201}\\
&=2\Bigr\{\sum_{n=-\i}^{\i}\df{q^{3n}}{1+q^{10n}}+\sum_{n=-\i}^{\i}\df{q^{5n+1}}{1+q^{10n+2}}\Bigl\}\label{d202}\\
&=1+\sum_{n=1}^{\i}\Bigr(\df{-20}{n}\Bigl)\df{q^n}{1-q^n}+\sum_{n=1}^{\i}\Bigr(\df{n}{5}\Bigl)\df{q^n}{1+q^{2n}}.\label{d203}
\end{align}

Furthermore,
\begin{equation}\label{d20p1}
1+\sum_{n=1}^{\i}\Bigr(\df{-20}{n}\Bigl)\df{q^n}{1-q^n}=\df{E(q^2)E(q^4)E(q^5)E(q^{10})}{
 E(q)E(q^{20}) }
\end{equation}
and
 \begin{equation}\label{d20p2}
\sum_{n=1}^{\i}\Bigr(\df{n}{5}\Bigl)\df{q^n}{1+q^{2n}}=q\df{E(q)E(q^2)E(q^{10})E(q^{20})}{
E(q^4)E(q^5) }.
\end{equation}

\end{theorem}
\begin{proof}
Employing \eqref{1psi1} with $q, a$ and $b$ replaced by
$q^{10}, q$ and $-1$, respectively, we find that
\begin{equation}
 \sum_{n=-\i}^{\i}\df{q^n}{1+q^{10n}}=E^3(q^{10})\df{f(q,q^9)}{f(-q,-q^9)f(1,q^{10})}.\label{l1}\\
\end{equation}
From \eqref{1psi1}, we similarly find that
\begin{equation}
\sum_{n=-\i}^{\i}\df{q^{5n+2}}{1+q^{10n+4}}=q^2E^3(q^{10})\df{f(q,q^9)}{f(-q^5,-q^5)f(q^4,q^6)}.\label{l2}
\end{equation}
From  \eqref{ff}, with $q$ replaced by $q^5$ and with $a=b=q^2$,
we obtain
\begin{equation}
f(-q^2,q^3)f(-q^2,q^3)=f(-q^5,-q^5)f(q^4,q^6)-q^2f(1,q^{10})f(-q,-q^9).\label{ff1}
\end{equation}
By \eqref{l1}, \eqref{l2} and \eqref{ff1}, we conclude that
\begin{align}
&\sum_{n=-\i}^{\i}\df{q^n}{1+q^{10n}}-\sum_{n=-\i}^{\i}\df{q^{5n+2}}{1+q^{10n+4}}\notag\\
&\;\;=E^3(q^{10})\df{f(q,q^9)}{f(-q,-q^9)f(1,q^{10})}-q^2E^3(q^{10})\df{f(q,q^9)}{f(-q^5,-q^5)f(q^4,q^6)}\notag\\
&\;\;=\df{E^3(q^{10})f(q,q^9)}{f(-q,-q^9)f(1,q^{10})f(-q^5,-q^5)f(q^4,q^6)}\Bigr\{f(-q^5,-q^5)f(q^4,q^6)-q^2f(1,q^{10})f(-q,-q^9)\Bigl\}\notag\\
&\;\;=\df{E^3(q^{10})f(q,q^9)}{f(-q,-q^9)f(1,q^{10})f(-q^5,-q^5)f(q^4,q^6)}f(-q^2,q^3)f(-q^2,q^3)\notag\\
&\;\;=\df{1}{2}\t(q)\t(q^5),\notag
\end{align}
after several applications of \eqref{19III}. Similarly, we
find that
\begin{align}
&\sum_{n=-\i}^{\i}\df{q^{3n}}{1+q^{10n}}+\sum_{n=-\i}^{\i}\df{q^{5n+1}}{1+q^{10n+2}}\notag\\
&\;\;=E^3(q^{10})\df{f(q^3,q^7)}{f(-q^3,-q^7)f(1,q^{10})}+qE^3(q^{10})\df{f(q^3,q^7)}{f(-q^5,-q^5)f(q^2,q^8)}\notag\\
&\;\;=\df{E^3(q^{10})f(q^3,q^7)}{f(-q^3,-q^7)f(1,q^{10})f(-q^5,-q^5)f(q^2,q^8)}\Bigr\{f(-q^5,-q^5)f(q^2,q^8)+qf(1,q^{10})f(-q^3,-q^7)\Bigl\}\notag\\
&\;\;=\df{E^3(q^{10})f(q^3,q^7)}{f(-q^3,-q^7)f(1,q^{10})f(-q^5,-q^5)f(q^2,q^8)}f(-q,q^4)f(-q,q^4) \notag \\
&\;\;=\df{1}{2}\t(q)\t(q^5).\notag
\end{align}

Before we move on we would like to make the following
\begin{remark}
The following generalized Lambert series identity for
$\t(q)\t(q^5)$ is given in \cite[cor.~6.5]{chan}

\begin{equation*}
\t(-q)\t(-q^5)=2\sum_{k=-\i}^{\i}\df{q^{k(5k+3)/2}}{1+q^{5k}}-2q\sum_{k=-\i}^{\i}\df{q^{k(5k+7)/2}}{1+q^{5k+2}}.
\end{equation*}
It would be interesting to find a direct proof that
\begin{equation*}
\sum_{k=-\i}^{\i}(-1)^k\Bigr(\df{q^k}{1+q^{10k}}-\df{q^{5k+2}}{1+q^{10k+4}}\Bigl)
=\sum_{k=-\i}^{\i}\df{q^{k(5k+3)/2}}{1+q^{5k}}-q\sum_{k=-\i}^{\i}\df{q^{k(5k+7)/2}}{1+q^{5k+2}}.
\end{equation*}
\end{remark}

Next, we prove \eqref{d203}.
\begin{align}
&\sum_{n=-\i}^{\i}\Bigr\{\df{q^{5n+1}}{1+q^{10n+2}}-\df{q^{5n+2}}{1+q^{10n+4}}\Bigl\}\notag\\
&=\sum_{n=0}^{\i}\Bigr\{\df{q^{5n+1}}{1+q^{10n+2}}-\df{q^{5n+2}}{1+q^{10n+4}}-\df{q^{5n+3}}{1+q^{10n+6}}+\df{q^{5n+4}}{1+q^{10n+8}}\Bigl\}\notag\\
&=\sum_{n=1}^{\i}\Bigr(\df{n}{5}\Bigl)\df{q^n}{1+q^{2n}}.\label{2.3.h}
\end{align}

Also,
\begin{align}
&\sum_{n=-\i}^{\i}\df{q^n+q^{3n}}{1+q^{10n}}\notag\\
&=1+\sum_{j \in
\{1,3,7,9\}}\sum_{n=1}^{\i}\df{q^{jn}}{1+q^{10n}}\notag\\
&=1+\sum_{j \in
\{1,3,7,9\}}\sum_{n=1}^{\i}\sum_{m=0}^{\i}(-1)^m q^{jn}q^{10nm} \notag\\
&=1+\sum_{j \in
\{1,3,7,9\}}\sum_{n=1}^{\i}\sum_{m=0}^{\i}(-1)^m q^{n(10m+j)} \notag\\
&=1+\sum_{j \in
\{1,3,7,9\}}\sum_{m=0}^{\i}(-1)^m \df{q^{10m+j}}{1-q^{10m+j}} \notag\\
&=1+\sum_{n=1}^{\i}\Bigr(\df{-20}{n}\Bigl)\df{q^n}{1-q^n}.\label{2.3.hh}
\end{align}
Now using \eqref{2.3.h} and \eqref{2.3.hh} together with
\eqref{d201} and \eqref{d202}, we see that \eqref{d203} is proved.

The equations \eqref{d20p1} and \eqref{d20p2}  are essentially
given in \cite[eqs.~3.2, 3.29]{Li}. Moreover, the $eta$-quotients
that appear in these equations are included in a list of certain
multiplicative functions determined by Y.~Martin \cite{martin}. We
should emphasize that \eqref{d201}--\eqref{d203} are new. We
observe directly that the coefficients of the two Lambert series
in \eqref{d203} are multiplicative and that they differ at most by
a sign. This enables us to compute the coefficients of
$\t(q)\t(q^5)$.

\end{proof}
\begin{corollary}
Let $a(n)$ be the number of representations of a positive integer
$n$ by  quadratic form $k^2+5l^2$. If the prime factorization of
$n$ is given by
\begin{equation*}
n=2^a5^b\prod_{i=1}^{r}p_i^{v_i}\prod_{j=1}^sq_j^{w_j},
\end{equation*}
where  $p_i \equiv 1, 3, 7$, \text{or} $ 9$ $\pmod {20}$ and $q_i
\equiv 11, 13, 17$, \text{or} $19$ $\pmod {20}$, then
\begin{equation}
a(n)=\Bigr( 1+(- 1)^{a + t  } \Bigl) \prod_{i=1}^{r} (1+v_i)
\prod_{j=1}^s\df{1+ (-1)^{w_j}}{2},\label{ffc}
\end{equation}
where $t$ is the number of prime factors of $n$, counting
multiplicity, that are congruent to  $3$ or $7$ $\pmod{20}$.
\end{corollary}
\begin{proof}
Observe that
\begin{align}
\sum_{n=1}^{\i}\Bigr(\df{n}{5}\Bigl)\df{q^n}{1+q^{2n}}&=\sum_{n=1}^{\i}\sum_{m=0}^{\i}(-1)^m\Bigr(\df{n}{5}\Bigl)q^{n(2m+1)}
=\sum_{n=1}^{\i}\sum_{m=1}^{\i}\Bigr(\df{-4}{m}\Bigl)\Bigr(\df{n}{5}\Bigl)q^{nm}\notag\\
&=\sum_{n=1}^{\i}\Bigr(\sum_{d|n}\Bigr(\df{-4}{d}\Bigl)\Bigr(\df{n/d}{5}\Bigl)\Bigl)q^n.
\end{align}
Similarly,
\begin{equation}
\sum_{n=1}^{\i}\Bigr(\df{-20}{n}\Bigl)\df{q^n}{1-q^n}=\sum_{n=1}^{\i}\Bigr(\sum_{d|n}\bigr(\df{-20}{d}\bigl)\Bigl)q^n.
\end{equation}
Define
\begin{equation}
b(n)=\sum_{d|n}\bigr(\df{-20}{d}\bigl) \;\; \text{and}\;\;
c(n)=\sum_{d|n}\Bigr(\df{-4}{d}\Bigl)\Bigr(\df{n/d}{5}\Bigl).
\end{equation}
We have by \eqref{d203} that $a(n)=b(n)+c(n)$. Clearly both $b(n)$
and $c(n)$ are multiplicative functions. Therefore, one only needs
to find their values at prime powers. It is easy to check that for
a prime $p$

\begin{equation}
 b(p^{\a})=\left\{ \begin{array}{ll}
         1  &  \text{if}\;\; p=2 \;\text{or}\; 5,\\
         1+\a &\text{if}\;\;p
\equiv 1, 3, 7, 9 \pmod{20},\\
\df{1+(-1)^{\a}}{2} & \text{if} \;\;p \equiv 11, 13, 17, 19
\pmod{20}, \end{array} \right. \label{bdet}
\end{equation}

and

\begin{equation}
 c(p^{\a})=\left\{ \begin{array}{ll}
         (-1)^{\a}  &  \text{if}\;\; p=2, \\
         1  &  \text{if}\;\; p=5, \\
         1+\a &\text{if}\;\;p
\equiv 1, 9 \pmod{20},\\
 (-1)^{\a}(1+\a) &\text{if}\;\;p
\equiv 3, 7 \pmod{20},\\
\df{1+(-1)^{\a}}{2} & \text{if} \;\;p \equiv 11, 13, 17, 19
\pmod{20}. \end{array} \right. \label{cdet}
\end{equation}

Equivalent reformulations of \eqref{ffc} can also be found in \cite[p.~84,
ex.~1]{D} and \cite[Thr.~7]{NH}. We should remark  that
\eqref{ffc} implies conjectures of Fermat and Euler for $(1,0,5)$
stated in the introduction. The last two equations immediately imply \eqref{ffc}.

\end{proof}

We now determine the representations of integers by the quadratic
form $(2,2,3)$ and make some further observations.

\begin{corollary}
Let $d(n)$ be a number of representations of a positive integer
$n$ by  the quadratic form $2k^2+2kl+3l^2$. If the prime
factorization of $n$ is given by
\begin{equation*}
n=2^a5^b\prod_{i=1}^{r}p_i^{v_i}\prod_{j=1}^sq_j^{w_j},
\end{equation*}
where  $p_i \equiv 1, 3, 7$, \text{or} $ 9$ $\pmod {20}$ and $q_i
\equiv 11, 13, 17$, \text{or} $19$ $\pmod {20}$, then
\begin{equation}
d(n)=\Bigr( 1-(- 1)^{a + t  } \Bigl) \prod_{i=1}^{r} (1+v_i)
\prod_{j=1}^s\df{1+ (-1)^{w_j}}{2},\label{ff2c}
\end{equation}
where $t$ is the number of prime factors of n, counting
multiplicity, that are congruent to  $3$ or $7\pmod{20}$.
\end{corollary}

\begin{proof}
Recall that
\begin{equation*}
\sum_{n=0}^{\i}a(n)q^n:=\sum_{k,\,l=-\i}^{\i}q^{k^2+5l^2}.
\end{equation*}
By comparing  \eqref{ffc} and \eqref{ff2c},  it suffices to show
that $d(n)=a(2n)$ for all $n$ $\in$ $\mathbb{N}$. To that end we
observe

\begin{align*}
\sum_{n=0}^{\i}d(n)q^n&=\sum_{n,m=-\i}^{\i}q^{2n^2+2nm+3m^2}\\
&=\sum_{n,m=-\i}^{\i}q^{2n^2+2n(2m)+3(2m)^2}+\sum_{n,m=-\i}^{\i}q^{2n^2+2n(2m+1)+3(2m+1)^2}\\
&=\sum_{n,m=-\i}^{\i}q^{2\bigr((n+m)^2+5m^2\bigl)}+\sum_{n,m=-\i}^{\i}q^{\tfrac{(2n+2m+1)^2+5(2m+1)^2}{2}}\\
&=\sum_{n,m=-\i}^{\i}q^{\frac{(2n)^2+5(2m)^2}{2}}+\sum_{n,m=-\i}^{\i}q^{\tfrac{(2n+1)^2+5(2m+1)^2}{2}}\\
&=\sum_{n=0}^{\i}a(2n)q^n.
\end{align*}
\end{proof}

By \eqref{d203}, \eqref{ffc} and \eqref{ff2c}  we find that
\begin{align}
\sum_{n,m=-\i}^{\i}q^{2n^2+2nm+3m^2}
=1+\sum_{n=1}^{\i}\Bigr(\df{-20}{n}\Bigl)\df{q^n}{1-q^n}-\sum_{n=1}^{\i}\Bigr(\df{n}{5}\Bigl)\df{q^n}{1+q^{2n}}.\label{d2032}
\end{align}
Also by adding identities in \eqref{d203} and \eqref{d2032},
we conclude that
\begin{equation*}
\sum_{n,m=-\i}^{\i}q^{n^2+5m^2}+\sum_{n,m=-\i}^{\i}q^{2n^2+2nm+3m^2}
=2+2\sum_{n=1}^{\i}\Bigr(\df{-20}{n}\Bigl)\df{q^n}{1-q^n}.
\end{equation*}
This last equation is a special case of Dirichlet's formula
\cite[p.~123, thr.~4]{Li2}. Comparing \eqref{ffc} and \eqref{ff2c}
we see that $a(n)d(n) =0.$ This means that a positive integer
cannot be represented by (1,0,5) and (2,2,3) at the same time.

We end this section by proving a Lambert series representation for
$\psi(q)\psi(q^5)$.
\begin{theorem}
\begin{equation}\label{sss}
\psi(q)\psi(q^5)=\sum_{n=-\i}^{\i}\df{q^{3n}+q^{7n+1}}{1-q^{20n+5}}=\sum_{n=-\i}^{\i}\df{q^{n}+q^{9n+6}}{1-q^{20n+15}}.
\end{equation}
\end{theorem}
\begin{proof}
By two applications of \eqref{1psi1} with $q$ replaced by
$q^{20}$, $a=q^3,\,q^7$ and $b=q^5$, we find that
\begin{align}
\sum_{n=-\i}^{\i}\df{q^{3n}+q^{7n+1}}{1-q^{20n+5}}&=E^3(q^{20})\df{f(-q^8,-q^{12})}{f(-q^3,-q^{17})f(-q^5,-q^{15})}
+qE^3(q^{20})\df{f(-q^8,-q^{12})}{f(-q^7,-q^{13})f(-q^5,-q^{15})}\notag\\
&=\df{E^3(q^{20})f(-q^8,-q^{12})}{f(-q^5,-q^{15})f(-q^3,-q^{17})f(-q^7,-q^{13})}
\Bigr\{f(-q^7,-q^{13})+qf(-q^3,-q^{17})\Bigl\}\notag\\
 &=\df{E^3(q^{20})f(-q^8,-q^{12})}{f(-q^5,-q^{15})f(-q^3,-q^{17})f(-q^7,-q^{13})}
f(q,-q^4),\label{nds}
\end{align}
where in the last step we use \eqref{2dd} with $a=q$ and $b=-q^4$.
It is now easy to verify by several applications of \eqref{19III}
that \eqref{nds} is equal to $\psi(q)\psi(q^5)$. The proof of the
second representation given in \eqref{sss} is very similar to that
of the first one and so we forego its proof.
\end{proof}

\bigskip

\section{Lambert Series Identities for $\sum_{n,m=-\i}^{\i}q^{n^2+6m^2}$  and $\sum_{n,m=-\i}^{\i}q^{2n^2+3m^2}$}

\begin{theorem}\label{DD24}
Let
\begin{equation}\label{p1p2def}
P(q):=\df{E(q^2)E(q^3)E(q^8)E(q^{12})}{E(q)E(q^{24})}
\;\;\text{and}\;\;
Q(q):=q\df{E(q)E(q^4)E(q^6)E(q^{24})}{E(q^3)E(q^8)}.
\end{equation}
Then,
\begin{align}
P(q)&=\sum_{n=-\i}^{\i}\df{q^n+q^{5n}}{1+q^{12n}}=1+\sum_{n=1}^{\i}\Bigr(\df{-6}{n}\Bigl)\df{q^n}{1-q^n},\label{forP}\\
Q(q)&=\sum_{n=-\i}^{\i}\df{q^{3n+1}-q^{9n+3}}{1+q^{12n+4}}=\sum_{n=1}^{\i}\Bigr(\df{n}{3}\Bigl)\df{q^n(1-q^{2n})}{1+q^{4n}}.\label{forQ}
\end{align}
Moreover,
\begin{align}
\t(q)\t(q^6)&=P(q)+Q(q)\label{d24}\\&=\sum_{n=-\i}^{\i}\df{q^n+q^{5n}}{1+q^{12n}}+\sum_{n=-\i}^{\i}\df{q^{3n+1}-q^{9n+3}}{1+q^{12n+4}}\label{d241}\\
&=2\Bigr\{\sum_{n=-\i}^{\i}\df{q^n}{1+q^{12n}}-\sum_{n=-\i}^{\i}\df{q^{9n+3}}{1+q^{12n+4}}\Bigl\}\label{16half}\\
&=2\Bigr\{\sum_{n=-\i}^{\i}\df{q^{5n}}{1+q^{12n}}+\sum_{n=-\i}^{\i}\df{q^{3n+1}}{1+q^{12n+4}}\Bigl\}\label{16ohalf}\\
&=1+\sum_{n=1}^{\i}\Bigr(\df{-6}{n}\Bigl)\df{q^n}{1-q^n}+\sum_{n=1}^{\i}\Bigr(\df{n}{3}\Bigl)\df{q^n(1-q^{2n})}{1+q^{4n}}\label{d242}\\
\intertext{and}\notag\\
\t(q^2)\t(q^3)&=P(q)-Q(q)\label{d24ii}\\&=\sum_{n=-\i}^{\i}\df{q^n+q^{5n}}{1+q^{12n}}-\sum_{n=-\i}^{\i}\df{q^{3n+1}-q^{9n+3}}{1+q^{12n+4}}\label{d243}\\
&=2\Bigr\{\sum_{n=-\i}^{\i}\df{q^n}{1+q^{12n}}-\sum_{n=-\i}^{\i}\df{q^{3n+1}}{1+q^{12n+4}}\Bigl\}\label{23half}\\
&=2\Bigr\{\sum_{n=-\i}^{\i}\df{q^{5n}}{1+q^{12n}}+\sum_{n=-\i}^{\i}\df{q^{9n+3}}{1+q^{12n+4}}\Bigl\}\label{23ohalf}\\
&=1+\sum_{n=1}^{\i}\Bigr(\df{-6}{n}\Bigl)\df{q^n}{1-q^n}-\sum_{n=1}^{\i}\Bigr(\df{n}{3}\Bigl)\df{q^n(1-q^{2n})}{1+q^{4n}}.\label{d244}
\end{align}
\end{theorem}

\begin{proof}
By employing \eqref{ff} with $a=q,\,b=q^{11},\,c=-q^5$ and
$d=-q^7$, we find that
\begin{align}
f(q,q^{11})f(-q^5,-q^7)&=f(-q^6,-q^{18})f(-q^8,-q^{16})+qf(-q^4,-q^{20})f(-q^6,-q^{18})\notag\\&=
\psi(-q^6)E(q^8)+q\psi(-q^6)f(-q^4,-q^{20}).\label{Ppre}
\end{align}

By two applications of \eqref{1psi1} with $q$ replaced by
$q^{12}$, $a=q, q^5$, $b=-1$, and by \eqref{Ppre}, we find that
\begin{align}
\sum_{n=-\i}^{\i}\df{q^n+q^{5n}}{1+q^{12n}}
&=E^3(q^{12})\df{f(q,q^{11})}{f(1,q^{12})f(-q,-q^{11})}+E^3(q^{12})\df{f(q^5,q^7)}{f(1,q^{12})f(-q^5,-q^7)}\notag\\
&=\df{E^3(q^{12})}{f(1,q^{12})f(-q,-q^{11})f(-q^5,-q^7)}\Bigr\{f(q,q^{11})f(-q^5,-q^7)+f(-q,-q^{11})f(q^5,q^7)\Bigl\}\notag\\
&=2\df{E^3(q^{12})}{f(1,q^{12})f(-q,-q^{11})f(-q^5,-q^7)}\psi(-q^6)E(q^8)\notag\\
&=\df{E(q^2)E(q^3)E(q^8)E(q^{12})}{E(q)E(q^{24})},\label{need1}
\end{align}
after several applications of \eqref{19III}.

By \eqref{2dd}, we observe that
\begin{equation}\label{edd}
E(q)=f(-q,-q^2)=f(q^5,q^7)-qf(q,q^{11}).
\end{equation}

Arguing as above and using \eqref{edd}, we conclude that
\begin{align}
\sum_{n=-\i}^{\i}\df{q^{3n+1}-q^{9n+3}}{1+q^{12n+4}}
&=qE^3(q^{12})\df{f(q^5,q^{7})}{f(-q^3,-q^{9})f(q^4,q^{8})}-q^3E^3(q^{12})\df{f(q^{-1},q^{13})}{f(-q^3,-q^{9})f(q^4,q^8)}\notag\\
&=qE^3(q^{12})\df{f(q^5,q^{7})}{f(-q^3,-q^{9})f(q^4,q^{8})}-q^2E^3(q^{12})\df{f(q,q^{11})}{f(-q^3,-q^{9})f(q^4,q^8)}\notag\\
&=q\df{E^3(q^{12})}{\psi(-q^3)f(q^4,q^8)}\Bigr\{f(q^5,q^7)-qf(q,q^{11})\Bigl\}\notag\\
&=q\df{E^3(q^{12})E(q)}{\psi(-q^3)f(q^4,q^8)}
=q\df{E(q)E(q^4)E(q^6)E(q^{24})}{E(q^3)E(q^8)},\label{need2}
\end{align}
after several applications of \eqref{19III}.

The proofs of second part of \eqref{forP} and that of \eqref{forQ}
are similar to those of \eqref{2.3.h}, \eqref{2.3.hh} and so we omit
their proofs.

We prove \eqref{d24} and \eqref{d24ii} simultaneously by proving
\begin{align}
2P&=\t(q)\t(q^6)+\t(q^2)\t(q^3) \label{ch1}\\
\intertext{and}\notag\\
 2Q&=\t(q)\t(q^6)-\t(q^2)\t(q^3).\label{ch2}
\end{align}

First we prove \eqref{ch2}. We will need two identities of
Ramanujan \cite[p.~232]{III}, namely
\begin{equation}
2\df{\psi^3(q)}{\psi(q^3)}=\df{\t^3(q)}{\t(q^3)}+\df{\t^3(-q^2)}{\t(-q^6)}\label{ntr1}
\end{equation}
and
\begin{equation}
4q\psi(q^2)\psi(q^6)=\t(q)\t(q^3)-\t(-q)\t(-q^3).\label{trii}
\end{equation}
From \eqref{trii} with \eqref{tdd}, we find that
\begin{align}
4q\psi(q^2)\psi(q^6)&=\t(q)\t(q^3)-\t(-q)\t(-q^3)\notag\\
&=\bigr(\t(q^{4})+q\psi(q^{8})\bigl)\bigr(\t(q^{12})+
q^3\psi(q^{24})\bigl)-\bigr(\t(q^{4})-q\psi(q^{8})\bigl)\bigr(\t(q^{12})-q^3\psi(q^{24})\bigl)\notag\\
&=4q\Bigr\{\psi(q^8)\t(q^{12})+q^2\t(q^4)\psi(q^{24})\Bigl\}.\label{dssg}
\end{align}
Upon replacing $q^2$ by $q$ in \eqref{dssg}, we conclude that
\begin{equation}
\psi(q)\psi(q^3)=\psi(q^4)\t(q^{6})+q\t(q^2)\psi(q^{12}).\label{yt1}
\end{equation}
Similarly,
\begin{align}
&\t(q)\t(-q^3)-\t(-q)\t(q^3)\notag\\
&=\bigr(\t(q^{4})+q\psi(q^{8})\bigl)\bigr(\t(q^{12})-
q^3\psi(q^{24})\bigl)-\bigr(\t(q^{4})-q\psi(q^{8})\bigl)\bigr(\t(q^{12})+q^3\psi(q^{24})\bigl)\notag\\
&=4q\Bigr\{\psi(q^8)\t(q^{12})-q^2\t(q^4)\psi(q^{24})\Bigl\}\notag\\
&=4q\psi(-q^2)\psi(-q^6),\label{yt2}
\end{align}
where in the last step we used \eqref{yt1} with $q$ replaced by
$-q^2$. We are now ready to prove \eqref{ch2}. Recall that $Q(q)$
is defined by \eqref{p1p2def}. By several applications of
\eqref{19III}, we see that \eqref{ch2} is equivalent to
\begin{equation}
2q\df{\psi(-q)\psi(-q^2)\psi(-q^3)\psi(-q^6)}{\psi(q^4)\t(-q^3)}=\t(q)\t(q^6)-\t(q^2)\t(q^3),
\end{equation}
or
\begin{align}
2q\psi(-q)\psi(-q^2)\psi(-q^3)\psi(-q^6)&=\t(q)\t(q^6)\psi(q^4)\t(-q^3)-\t(q^2)\t(q^3)\psi(q^4)\t(-q^3)\notag\\
&=\t(q)\t(-q^3)\psi(q^4)\t(q^6)-\psi^2(q^2)\t^2(-q^6),\label{tgggg}
\end{align}
where we used the trivial identities
\begin{equation}
\t(q)\t(-q)=\t^2(-q^2)\;\; \text{and}
\;\;\psi^2(q)=\psi(q^2)\phi(q).
\end{equation}

When we employ \eqref{yt2} on the far left hand side of
\eqref{tgggg} and \eqref{yt1} on the right hand side of
\eqref{tgggg}, we find that

\begin{equation}
\psi(-q)\psi(-q^3)\bigr(\t(q)\t(-q^3)-\t(-q)\t(q^3)\bigl)=\bigr(\psi(q)\psi(q^3)+\psi(-q)\psi(-q^3)\bigl)\t(q)\t(-q^3)-2\psi^2(q^2)\t^2(-q^6).
\end{equation}
Upon cancellation, we see that
\begin{equation}
-\psi(-q)\psi(-q^3)\t(-q)\t(q^3)=\psi(q)\psi(q^3)\t(q)\t(-q^3)-2\psi^2(q^2)\t^2(-q^6).\label{almost}
\end{equation}
Next we multiply both sides of \eqref{almost} with
$\df{\t^2(q)}{\psi(q)\psi(q^3)\t^2(-q^6)}$ and obtain after
several applications of \eqref{19III} that
\begin{equation}
-\df{\t^3(-q^2)}{\t(-q^6)}=\df{\t^3(q)}{\t(q^3)}-2\df{\psi^3(q)}{\psi(q^3)},
\end{equation}
which is \eqref{ntr1}. Hence the proof of \eqref{ch2} is complete.

The proof of  \eqref{ch1} is very similar to that of \eqref{ch2}.
Recall that $Q(q)$ is defined by \eqref{p1p2def}. By several
applications of \eqref{19III}, we see that \eqref{ch1} is
equivalent to
\begin{equation}
2\df{\psi(-q)\psi(-q^2)\psi(-q^3)\psi(-q^6)}{\psi(q^{12})\t(-q)}=\t(q)\t(q^6)+\t(q^2)\t(q^3),
\end{equation}
or
\begin{align}
\psi(-q)\psi(-q^2)\psi(-q^3)\psi(-q^6)&=\t(q)\t(q^6)\psi(q^{12})\t(-q)+\t(q^2)\t(q^3)\psi(q^{12})\t(-q)\notag\\
&=\t^2(-q^2)\psi^2(q^6)+\t(-q)\t(q^3)\t(q^2)\psi(q^{12}).\label{tgggg2}
\end{align}

If we employ \eqref{yt2} on the far left hand side of
\eqref{tgggg2}, and \eqref{yt1} on the right hand side of
\eqref{tgggg2} and multiply both sides by $2q$, we find that

\begin{equation}
\psi(-q)\psi(-q^3)\bigr(\t(q)\t(-q^3)-\t(-q)\t(q^3)\bigl)=2q\psi^2(q^{6})\t^2(-q^2)+\t(-q)\t(q^3)\bigr(\psi(q)\psi(q^3)-\psi(-q)\psi(-q^3)\bigl).
\end{equation}
Upon cancellation, we find that
\begin{equation}
\psi(-q)\psi(-q^3)\t(q)\t(-q^3)=\psi(q)\psi(q^3)\t(-q)\t(q^3)+2q\psi^2(q^{6})\t^2(-q^2).\label{almost2}
\end{equation}
It is easy to see that \eqref{almost2} and \eqref{almost} are
``reciprocals'' of each other. For related definitions and modular
equations corresponding to \eqref{ch1} and \eqref{ch2} see
\cite[p.~230, Entry 5 (i)]{III}. Hence, the proof of \eqref{ch1}
is complete.

As an immediate corollary of  \eqref{ch1} and \eqref{ch2}, we note
the following two interesting theta function identities
\begin{equation}
\df{\t(q)\t(q^6)-\t(q^2)\t(q^3)}{\t(q)\t(q^6)+\t(q^2)\t(q^3)}=q\df{\t(-q)\psi(q^{12})}{\t(-q^3)\psi(q^4)}\label{yhhh}
\end{equation}
and
\begin{equation*}
\t^2(q)\t^2(q^6)-\t^2(q^2)\t^2(q^3)=4qE(q^2)E(q^4)E(q^6)E(q^{12})=4q\psi(q)\psi(-q)\psi(-q^{3})\psi(-q^{6}).
\end{equation*}

The identities \eqref{forP} and \eqref{forQ} together with
\eqref{d24} and \eqref{d24ii} clearly imply \eqref{d201} and
\eqref{d203}. To prove the remaining identities \eqref{16half},
\eqref{16ohalf}, \eqref{23half} and \eqref{23ohalf}, one only
needs to prove that
\begin{equation}\label{gghh}
\sum_{n=-\i}^{\i}\df{q^n-q^{5n}}{1+q^{12n}}=\sum_{n=-\i}^{\i}\df{q^{3n+1}+q^{9n+3}}{1+q^{12n+4}}.
\end{equation}
Arguing as in \eqref{need1} and \eqref{need2}, one can easily show
that
\begin{equation}
\sum_{n=-\i}^{\i}\df{q^n-q^{5n}}{1+q^{12n}}=\sum_{n=-\i}^{\i}\df{q^{3n+1}+q^{9n+3}}{1+q^{12n+4}}=\df{E(q^2)E(q^3)E(q^4)E(q^{24})}{E(q)E(q^8)}.
\end{equation}
Hence, the proof of the Theorem \ref{DD24} is complete.
\end{proof}

\begin{corollary}
Let $a(n)$ and $b(n)$ be the number of representations of a
positive integer $n$ by  quadratic form $k^2+6l^2$ and
$2k^2+3l^2$, respectively. If the prime factorization of $n$ is
given by
\begin{equation*}
n=2^a3^b\prod_{i=1}^{r}p_i^{v_i}\prod_{j=1}^sq_j^{w_j},
\end{equation*}
where  $p_i \equiv 1, 5, 7$, \text{or} $ 11$ $\pmod {24}$ and $q_i
\equiv 13, 17, 19$, \text{or} $23$ $\pmod {24}$, then
\begin{equation}
a(n)=\Bigr( 1+(- 1)^{a + b+t  } \Bigl) \prod_{i=1}^{r} (1+v_i)
\prod_{j=1}^s\df{1+ (-1)^{w_j}}{2}\label{24ffc}
\end{equation}
and
\begin{equation}
b(n)=\Bigr( 1-(- 1)^{a + b+t  } \Bigl) \prod_{i=1}^{r} (1+v_i)
\prod_{j=1}^s\df{1+ (-1)^{w_j}}{2},\label{24iiffc}
\end{equation}
where $t$ is a number of prime factors of $n$, counting
multiplicity, that are congruent to  $5$ or $11$ $\pmod{24}$.
\end{corollary}
\begin{proof}
Observe that
\begin{align*}
\sum_{n=1}^{\i}\Bigr(\df{n}{3}\Bigl)\df{q^n(1-q^{2n})}{1+q^{4n}}&=\sum_{n=1}^{\i}\sum_{m=0}^{\i}(-1)^m\Bigr(\df{n}{3}\Bigl)(q^{n(4m+1)}-q^{n(4m+3)})\notag\\
&=\sum_{n=1}^{\i}\sum_{m=1}^{\i}(-1)^m\Bigr(\df{n}{3}\Bigl)(q^{n(4m-1)}-q^{n(4m-3)})\notag\\
&=\sum_{n=1}^{\i}\sum_{m=1}^{\i}\Bigr(\df{m}{2}\Bigl)\Bigr(\df{n}{3}\Bigl)q^{nm}\notag\\
&=\sum_{n=1}^{\i}\Bigr(\sum_{d|n}\Bigr(\df{d}{2}\Bigl)\Bigr(\df{n/d}{3}\Bigl)\Bigl)q^n.
\end{align*}
 Similarly,
\begin{equation*}
\sum_{n=1}^{\i}\Bigr(\df{-6}{n}\Bigl)\df{q^n}{1-q^n}=\sum_{n=1}^{\i}\Bigr(\sum_{d|n}\bigr(\df{-6}{d}\bigl)\Bigl)q^n.
\end{equation*}
Define
\begin{equation*}
c(n)=\sum_{d|n}\bigr(\df{-24}{d}\bigl) \;\; \text{and}\;\;
d(n)=\sum_{d|n}\Bigr(\df{d}{2}\Bigl)\Bigr(\df{n/d}{3}\Bigl).
\end{equation*}
Equations \eqref{d242} and \eqref{d244} imply that $a(n)=c(n)+d(n)$ and
$b(n)=c(n)-d(n)$. Clearly both $c(n)$ and $d(n)$ are
multiplicative functions. Therefore, one only needs to find their
values at prime powers. It is easy to check that for a prime $p$

\begin{equation}
 c(p^{\a})=\left\{ \begin{array}{ll}
         1  &  \text{if}\;\; p=2 \;\text{or}\; 3,\\
         1+\a &\text{if}\;\;p
\equiv 1, 5, 7, 11 \pmod{24},\\
\df{1+(-1)^{\a}}{2} & \text{if} \;\;p \equiv 13, 17, 19, 23
\pmod{24} \end{array} \right.\label{24bdet}
\end{equation}

and

\begin{equation}
 d(p^{\a})=\left\{ \begin{array}{ll}
         (-1)^{\a}  &  \text{if}\;\; p=2 \;\text{or}\; 3, \\
         1+\a &\text{if}\;\;p
\equiv 1, 7 \pmod{24},\\
 (-1)^{\a}(1+\a) &\text{if}\;\;p
\equiv 5, 11 \pmod{24},\\
\df{1+(-1)^{\a}}{2} & \text{if} \;\;p \equiv 13, 17, 19, 23
\pmod{24}. \end{array} \right. \label{24cdet}
\end{equation}

From these two equations \eqref{24ffc} and \eqref{24iiffc} are
immediate. Equivalent reformulations of \eqref{24ffc} and
\eqref{24iiffc}  can also be found in \cite[p.~84, ex.~2]{D} and
\cite[Thr.~7]{NH}.
\end{proof}

\bigskip

\section{Lambert Series Identities for $\sum_{n,m=-\i}^{\i}q^{n^2+15m^2}$ and $\sum_{n,m=-\i}^{\i}q^{3n^2+5m^2}$}
\begin{theorem}
Let
\begin{equation}\label{60pqdef}
P(q):=\df{E(q)E(q^6)E(q^{10})E(q^{15})}{E(q^2)E(q^{30})}\;\;\text{and}\;\;Q(q):=q\df{E(q^2)E(q^3)E(q^{5})E(q^{30})}{E(q^6)E(q^{10})}.
\end{equation}
Then,
\begin{align}
P(q)&=1-\sum_{n=1}^{\i}\Bigr(\df{-15}{n}\Bigl)\df{q^n}{1+q^n},\label{60forP}\\
Q(q)&=\sum_{n=1}^{\i}\Bigr(\df{5}{n}\Bigl)\df{q^n(1+q^{n})}{1+q^{3n}}.\label{60forQ}
\end{align}
Moreover,
\begin{align}
\t(-q)\t(-q^{15})&=P(q)-Q(q)\label{d60}\\
&=1-\sum_{n=1}^{\i}\Bigr(\df{-15}{n}\Bigl)\df{q^n}{1+q^n}-\sum_{n=1}^{\i}\Bigr(\df{5}{n}\Bigl)\df{q^n(1+q^{n})}{1+q^{3n}}\label{d601}\\
\intertext{and}\notag\\
\t(-q^3)\t(-q^5)&=P(q)+Q(q)\label{d60ii}\\
&=1-\sum_{n=1}^{\i}\Bigr(\df{-15}{n}\Bigl)\df{q^n}{1+q^n}+\sum_{n=1}^{\i}\Bigr(\df{5}{n}\Bigl)\df{q^n(1+q^{n})}{1+q^{3n}}.\label{d602}
\end{align}
\end{theorem}
\begin{proof}
The identities \eqref{60forP}, \eqref{d60}, and \eqref{d60ii} were
observed by Ramanujan \cite[p.~379, Entry 10 (vi)]{III},
\cite[eq.~50]{wl} and \cite[p.~377, Entry 9 (v), (vi)]{III}. We
prove \eqref{60forQ}. \\
It is easy to observe that
\begin{equation}\label{60pre1}
\sum_{n=1}^{\i}\Bigr(\df{5}{n}\Bigl)\df{q^n(1+q^{n})}{1+q^{3n}}
=\sum_{n=-\i}^{\i}\df{q^{5n+1}+q^{10n+2}}{1+q^{15n+3}}-\sum_{n=-\i}^{\i}\df{q^{5n+1}+q^{10n+2}}{1+q^{15n+3}}.
\end{equation}
Next by four applications of \eqref{1psi1} on the right hand side
of \eqref{60pre1}, we have
\begin{align}
\sum_{n=1}^{\i}\Bigr(\df{5}{n}\Bigl)\df{q^n(1+q^{n})}{1+q^{3n}}&=E^3(q^{15})\Bigr\{
q\df{f(q^7,q^8)}{E(q^5)f(q^3,q^{12})}+q^2\df{f(q^2,q^{13})}{E(q^5)f(q^3,q^{12})}
\Bigl\}\notag\\
&\;\;\;-E^3(q^{15})\Bigr\{
q^2\df{f(q^4,q^{11})}{E(q^5)f(q^6,q^{9})}+q^3\df{f(q,q^{14})}{E(q^5)f(q^6,q^{9})}
\Bigl\}\notag\\
&=q\df{E^3(q^{15})}{E(q^5)f(q^3,q^{12})f(q^6,q^9)}\Bigr\{\bigr(f(q^7,q^8)+qf(q^2,q^{13})\bigl)f(q^6,q^9)\notag\\
&\;\;\;-q\bigr(f(q^4,q^{11})+qf(q,q^{14})\bigl)f(q^3,q^{12})
 \Bigl\}.\label{60pre2}
\end{align}
Now we employ \eqref{ff} for each term of \eqref{60pre2} inside
the parenthesis, the identity in \eqref{60pre2} now becomes
\begin{align}
&q\df{E^3(q^{15})}{E(q^5)f(q^3,q^{12})f(q^6,q^9)}\Bigr\{\bigr(f(q^{14},q^{16})-q^2f(q^{4},q^{26})\bigl)\bigr(f(q^{13},q^{17})-qf(q^7,q^{23})\bigl)\notag\\
&\;\;\;+q\bigr(f(q^8,q^{22})-q^2f(q^2,q^{28})\bigl)\bigr(f(q^{11},q^{19})-q^3f(q,q^{29})\bigl).\label{60pre3}
\end{align}
Recall that the Rogers-Ramanujan functions are defined by
\begin{equation}\label{GH}
G(q) :=\sum_{n=0}^{\i}\df{q^{n^2}}{(q;q)_n} \qquad \text{and}
\qquad H(q) :=\sum_{n=0}^{\i}\df{q^{n(n+1)}}{(q;q)_n}.
\end{equation}

These functions satisfy the famous Rogers--Ramanujan identities
\cite[pp.~214--215]{cp}
\begin{equation}\label{rridents}
G(q) = \df{1}{(q;q^5)_{\i}(q^4;q^5)_{\i}} \qquad \text{and} \qquad
H(q) = \df{1}{(q^2;q^5)_{\i}(q^3;q^5)_{\i}}.
\end{equation}
Our proof makes use of one of Ramanujan's forty identities for the
Rogers-Ramanujan functions, namely \cite{watson}
\begin{equation}
G(q)G(q^4)-qH(q)H(q^4)=\df{\t(q^5)}{E(q^2)}.\label{40gh}
\end{equation}

Next, we employ the quintuple product identity \eqref{qtp}, with
$q$ replaced by $q^{10}$ and $a=-q$ to find that
\begin{equation}
f(-q^{13},-q^{17})+qf(-q^7,-q^{23}) =
E(q^{10})\df{f(-q^2,-q^8)}{f(-q,-q^9)}=E(q^2)G(q).\label{qp1}
\end{equation}
Similarly, from \eqref{qtp} we find
\begin{align}
E(q^2)H(q)&=f(-q^{11},-q^{19})+q^3f(-q,-q^{29}),\label{qp2}\\
E(q)G(q^2)&=f(q^7,q^8)-qf(q^2,q^{13}),\label{qp3}\\
E(q)H(q^2)&=f(q^4,q^{11})-qf(q,q^{14}).\label{qp4}
\end{align}
Next, making use of \eqref{60pre1}-\eqref{60pre3}, \eqref{qp1}--\eqref{qp4},
\eqref{40gh}, and \eqref{60pqdef}, we conclude that

\begin{align}
&\sum_{n=1}^{\i}\Bigr(\df{5}{n}\Bigl)\df{q^n(1+q^{n})}{1+q^{3n}}\notag\\
&=q\df{E^3(q^{15})}{E(q^5)f(q^3,q^{12})f(q^6,q^9)}E^2(q^2)\Bigr\{G(q^4)G(-q)+qH(q^4)H(-q)\Bigl\}\notag\\
&=q\df{E^3(q^{15})}{E(q^5)f(q^3,q^{12})f(q^6,q^9)}E^2(q^2)\df{\t(-q^5)}{E(q^2)}\notag\\
&=q\df{E(q^2)E(q^3)E(q^{5})E(q^{30})}{E(q^6)E(q^{10})}\notag\\
&=Q(q).\notag
\end{align}
\end{proof}
Adding together \eqref{d60} and \eqref{d60ii} and replacing $q$ by
$-q$, we find that
\begin{equation*}
\t(q)\t(q^{15})+\t(q^3)\t(q^5)=2-2\sum_{n=1}^{\i}\Bigr(\df{-15}{n}\Bigl)\df{(-q)^n}{1+(-q)^n}.
\end{equation*}

It is instructive to compare this formula with  an equivalent
formula (50)  in \cite{wl} which states that
\begin{equation*}
\t(q)\t(q^{15}) +\t(q^3)\t(q^5)
=2+\sum_{n=1}^{\i}\tilde{a}(n)\df{q^n}{1-q^n},
\end{equation*}

where
\begin{equation*}
\tilde{a}(n) = 2\Bigr(\df{-60}{n}\Bigl)  -2\delta
(2|n)\Bigr(\df{-60}{n/2}\Bigl)
+2\delta(4|n)\Bigr(\df{-15}{n/4}\Bigl),
\end{equation*}
with
\begin{equation*}
\delta(a|b)=\left\{ \begin{array}{ll}
         1  &  \text{if} \;a|b,\\
         0 & \text{otherwise}\,.  \end{array} \right.
\end{equation*}

\begin{corollary}
Let $a(n)$ and $b(n)$ be the number of representations of a
positive integer $n$ by  quadratic form $k^2+15l^2$ and
$3k^2+5l^2$, respectively. If the prime factorization of $n$ is
given by
\begin{equation*}
n=2^a3^b5^c\prod_{i=1}^{r}p_i^{v_i}\prod_{j=1}^sq_j^{w_j},
\end{equation*}
where  $p_i \equiv 1, 2, 4$, \text{or} $ 8$ $\pmod {15}$, $p_i
\neq 2$ and $q_i \equiv 7, 11, 13$, \text{or} $14$ $\pmod {15}$,
then
\begin{equation}
a(n)=|a-1|\Bigr( 1+(- 1)^{a + b+c+t  } \Bigl) \prod_{i=1}^{r}
(1+v_i) \prod_{j=1}^s\df{1+ (-1)^{w_j}}{2}\label{60ffc}
\end{equation}
and
\begin{equation}
b(n)=|a-1|\Bigr( 1-(- 1)^{a + b+c+t  } \Bigl) \prod_{i=1}^{r}
(1+v_i) \prod_{j=1}^s\df{1+ (-1)^{w_j}}{2},\label{60iiffc}
\end{equation}
where $t$ is a number of odd prime factors of $n$, counting
multiplicity, that are congruent to  $2$ or $8$ $\pmod{15}$.
\end{corollary}
\begin{proof}
Observe that
\begin{align*}
Q(q)=\sum_{n=1}^{\i}\Bigr(\df{5}{n}\Bigl)\df{q^n(1+q^{n})}{1+q^{3n}}&
=\sum_{n=1}^{\i}\sum_{m=0}^{\i}(-1)^m\Bigr(\df{5}{n}\Bigl)(q^{n(3m+1)}+q^{n(3m+2)})\notag\\
&=-\sum_{n=1}^{\i}\sum_{m=1}^{\i}(-1)^m\Bigr(\df{5}{n}\Bigl)(q^{n(3m-1)}+q^{n(3m-2)})\notag\\
&=-\sum_{n=1}^{\i}\sum_{m=1}^{\i}(-1)^m\Bigr(\df{-3}{m}\Bigl)\Bigr(\df{5}{n}\Bigl)q^{nm}\notag\\
&=-\sum_{n=1}^{\i}\Bigr(\sum_{d|n}(-1)^d\Bigr(\df{-3}{d}\Bigl)\Bigr(\df{5}{n/d}\Bigl)\Bigl)q^n.
\end{align*}
Therefore,
\begin{equation}
Q(-q)=-\sum_{n=1}^{\i}\Bigr(\sum_{d|n}(-1)^{n+d}\Bigr(\df{-3}{d}\Bigl)\Bigr(\df{5}{n/d}\Bigl)\Bigl)q^n.\label{fix77a}
\end{equation}
 Similarly,
\begin{equation}
P(-q)=1+\sum_{n=1}^{\i}\Bigr(\sum_{d|n}(-1)^{n+d}\Bigr(\df{-15}{n/d}\Bigl)\Bigl)q^n.\label{fix77b}
\end{equation}
Now we define
\begin{equation*}
c(n)=\sum_{d|n}(-1)^{n+d}\Bigr(\df{-15}{n/d}\Bigl) \;\;
\text{and}\;\;
d(n)=\sum_{d|n}(-1)^{n+d}\Bigr(\df{-3}{d}\Bigl)\Bigr(\df{5}{n/d}\Bigl).
\end{equation*}
By \eqref{d60} and \eqref{d60ii} we have $a(n)=c(n)+d(n)$ and
$b(n)=c(n)-d(n)$, for $n>0$.  Using the fact that $(-1)^{n+1}$ is
a multiplicative function of $n$, we conclude that $c(n)$ and
$d(n)$ are multiplicative functions. From \eqref{60pqdef},
\eqref{fix77a} and \eqref{fix77b} we also observe that the
following $eta$-quotients are multiplicative
\begin{align}
P(-q)&=\df{E(-q)E(q^6)E(q^{10})E(-q^{15})}{E(q^2)E(q^{30})}=\df{E^2(q^2)E(q^6)E(q^{10})E^2(q^{30})}{E(q)E(q^4)E(q^{15})E(q^{60})},\label{fix60i}\\
-Q(-q)&=q\df{E(q^2)E(-q^3)E(-q^{5})E(q^{30})}{E(q^6)E(q^{10})}=q\df{E(q^2)E^2(q^6)E^2(q^{10})E(q^{30})}{E(q^3)E(q^5)E(q^{12})E(q^{20})},\label{fix60ii}\\
Q(q)&=q\df{E(q^2)E(q^3)E(q^{5})E(q^{30})}{E(q^6)E(q^{10})}.\label{fix60iii}
\end{align}

It is easy to check that for a prime $p$

\begin{equation}
 c(p^{\a})=\left\{ \begin{array}{ll}
         |\a-1|  &  \text{if}\;\; p=2, \\
         1  &  \text{if}\;\; p=3 \;\text{or}\; 5,\\
         1+\a &\text{if}\;\;p
\equiv 1, 2, 4, 8 \pmod{15}, p \neq 2, \\
\df{1+(-1)^{\a}}{2} & \text{if} \;\;p \equiv 7, 11, 13, 14
\pmod{15}, \end{array} \right. \label{64bdet}
\end{equation}
and
\begin{equation}
 d(p^{\a})=\left\{ \begin{array}{ll}
(-1)^{\a}|\a-1|  &  \text{if}\;\; p=2, \\
         (-1)^{\a}  &  \text{if}\;\; p=3 \;\text{or}\; 5, \\
         1+\a &\text{if}\;\;p
\equiv 1, 4 \pmod{15},\\
 (-1)^{\a}(1+\a) &\text{if}\;\;p
\equiv 2, 8 \pmod{15}, p \neq 2, \\
\df{1+(-1)^{\a}}{2} & \text{if} \;\;p \equiv 7, 11, 13, 14
\pmod{15}. \end{array} \right. \label{64cdet}
\end{equation}

From these two equations \eqref{60ffc} and \eqref{60iiffc} are
immediate. Equivalent reformulations of \eqref{60ffc} and
\eqref{60iiffc}  can also be found in \cite{NH}.
\end{proof}

\bigskip

\section{Representations by the the quadratic form $k^2+27l^2$}

In this section, we give a formula for the number of
representations of a positive integer by the quadratic form
$k^2+27l^2$.

\begin{theorem}
\begin{equation}
\t(q)\t(q^{27})=\df{\t(q)\t(q^3)-\t(q^3)\t(q^9)}{3}+\t(q^9)\t(q^{27})+\df{4}{3}qE(q^6)E(q^{18}).\label{tdsss}
\end{equation}
Let $a(n)$, $b(n)$ be the number of representations of a positive
integer $n$ by  quadratic forms $(1,0,27)$ and $(4,2,7)$,
respectively.\\
If $n \not \equiv 1 \pmod 6$, then
\begin{equation}\label{1not6case}
a(n)=b(n)=\left\{ \begin{array}{ll}
        \bigr( 3-2\delta_{\a,0} \bigl)(1+(-1)^{\a})\prod_{i=1}^{r} (1+v_i) \prod_{j=1}^s\df{1+
(-1)^{w_j}}{2} &  \text{if}\;\; \b \geq 2,\\
     (1+(-1)^{\a})\prod_{i=1}^{r} (1+v_i) \prod_{j=1}^s\df{1+
(-1)^{w_j}}{2}&\text{if} \;\; \b=0\;\; \text{and}\;\; \a >0,\\ 0 &
\text{otherwise},\end{array} \right.
\end{equation}
where $n$ has the prime factorization
\begin{equation*}
n=2^{\a}3^{\b}\prod_{i=1}^{r}p_i^{v_i}\prod_{j=1}^sq_j^{w_j},
\end{equation*}
with  $p_i \equiv 1 \pmod {3}$ and $2\neq q_i \equiv 2 \pmod {3}$
and
\begin{equation}
\delta_{j,0}:=\left\{ \begin{array}{ll}
        1 &  \text{if}\;\; j=0,\\
         0 &\text{otherwise} .\end{array} \right.
\end{equation}

If $ n \equiv 1 \pmod 6$, then

\begin{align}
a(n)&=\df{2}{3}\prod_{i=1}^{r} (1+v_i)\Bigr(\prod_{i=1}^{s}
(1+u_i)+2\prod_{i=1}^{s} \Bigr((1+u_i)\mid 3\Bigl)\Bigl)
\prod_{i=1}^{t}\df{1+ (-1)^{w_i}}{2},\label{frman}\\
b(n)&=\df{2}{3}\prod_{i=1}^{r} (1+v_i)\Bigr(\prod_{i=1}^{s}
(1+u_i)-\prod_{i=1}^{s} \Bigr((1+u_i)\mid 3)\Bigl)\Bigl)
\prod_{i=1}^{t}\df{1+ (-1)^{w_i}}{2},\label{frmbn}
\end{align}
where $n$ has the prime factorization
\begin{equation}
\prod_{i=1}^{r}p_i^{v_i}\prod_{i=1}^{s}q_i^{u_i}\prod_{i=1}^{t}Q_j^{w_i},
\end{equation}
with  $p_i \equiv 1 \pmod {3}$, $2^{\tfrac{p_i-1}{3}} \equiv 1
\pmod {p_i}$,  $ q_i \equiv 1 \pmod {3}$, $2^{\tfrac{q_i-1}{3}}
\not \equiv 1 \pmod {q_i}$ and
 $ 2 \neq Q_i \equiv 2 \pmod {3}$.
\end{theorem}
\begin{proof}
Observe that
\begin{align}
&\sum_{u,v=-\i}^{\i}q^{7u^2+2uv+4v^2}=\sum_{k=0}^{3}\sum_{\;s,v=-\i}^{\i}q^{7(4s+k)^2+2(4s+k)v+4v^2}\notag\\
&=\sum_{k=0}^{3}\sum_{\;s,r=-\i}^{\i}q^{7(4s+k)^2+2(4s+k)(r-s)+4(r-s)^2}
=\sum_{k=0}^{3}q^{7k^2}\sum_{\;r=-\i}^{\i}q^{2(2r^2+kr)}\sum_{\;s=-\i}^{\i}q^{54(2s^2+ks)}\notag\\
&=f(q^4,q^4)f(q^{108},q^{108})+2q^7f(q^2,q^6)f(q^{54},q^{162})+q^{28}f(1,q^8)f(1,q^{216})\notag\\
&=(\t(q)\t(q^{27})+\t(-q)\t(-q^{27}))/2+2q^7\psi(q^2)\psi(q^{54}),\label{lazim1}
\end{align}
where in the last step we used \eqref{tdd}.

Similarly, we find that

\begin{align}
&\sum_{u,v=-\i}^{\i}q^{7u^2+2uv+4v^2}=\sum_{u,v=-\i}^{\i}q^{7(u-v)^2+2(u-v)v+4v^2}=\sum_{u,v=-\i}^{\i}q^{7u^2-12uv+9v^2}\notag\\
&=\sum_{k=0}^{2}\sum_{\;s,v=-\i}^{\i}q^{7(3s+k)^2-12(3s+k)v+9v^2}
=\sum_{k=0}^{2}\sum_{\;s,r=-\i}^{\i}q^{7(3s+k)^2-12(3s+k)(r+2s)+9(r+2s)^2}\notag\\
&=\sum_{k=0}^{2}q^{7k^2}\sum_{\;r=-\i}^{\i}q^{3(3r^2-4kr)}\sum_{\;s=-\i}^{\i}q^{9(3s^2+2ks)}\notag\\
&=f(q^9,q^9)f(q^{27},q^{27})+2q^4f(q^3,q^{15})f(q^9,q^{45})\notag\\
&=\t(q^9)\t(q^{27})+(\t(q)-\t(q^9))(\t(q^3)-\t(q^{27}))/2\notag\\
&=\bigr(3\t(q^9)\t(q^{27})+\t(q)\t(q^{3})-\t(q)\t(q^{27})-\t(q^3)\t(q^{9})\bigl)/2,\label{lazim4}
\end{align}
where we used \eqref{tdd3}.

Lastly, we  need the following identity of Ramanujan \cite[p.~359,
Entry.~4, (iv)]{III}

\begin{equation}
(\t(q)\t(q^{27})-\t(-q)\t(-q^{27}))/2-2q^7\psi(q^2)\psi(q^{54})=2qE(q^6)E(q^{18}).\label{lazim2}
\end{equation}

From \eqref{lazim1}--\eqref{lazim2}, we find that
\begin{equation}
\t(q)\t(q^{27})=\bigr(3\t(q^9)\t(q^{27})+\t(q)\t(q^{3})-\t(q)\t(q^{27})-\t(q^3)\t(q^{9})\bigl)/2+2qE(q^6)E(q^{18}),
\end{equation}
which is \eqref{tdsss}. The formulas
\eqref{lazim1}--\eqref{lazim2} are special cases of more general
formula  \cite[thr.~(3.1), cor.~ (3.3)]{RG}. In fact, Ramanujan's
identity, \eqref{lazim2}, can be stated as follows
\begin{align*}
&\sum_{u,v=\i}^{\i}(-1)^{u+v}q^{(7(2u+1)^2-12(2u+1)(2v+1)+9(2v+1)^2)/4}\\
&=(\t(q)\t(q^{27})-\t(-q)\t(-q^{27}))/2-2q^7\psi(q^2)\psi(q^{54})\\
&=2qE(q^6)E(q^{18}).
\end{align*}

Recall that
\begin{equation}
1+\sum_{n=1}^{\i}a(n)q^n:=\t(q)\t(q^{27}) \;\;\;\text{and}\;\;\;
1+\sum_{n=1}^{\i}b(n)q^n:=\sum_{n,m=-\i}^{\i}q^{4n^2+2nm+7m^2}.
\end{equation}
We also define $c(n)$ and $d(n)$ by
\begin{equation}
1+\sum_{n=1}^{\i}c(n)q^n=\t(q)\t(q^{3})
\;\;\;\text{and}\;\;\;\sum_{n=1}^{\i}d(n)q^n=qE(q^6)E(q^{18}).
\end{equation}
Using the following Lambert series expansion for $\t(q)\t(q^3)$
(see for example \cite[p.~75, eq.~(3.7.8)]{berndt2})
\begin{equation}
\t(q)\t(q^3)=1+2\sum_{n=1}^{\i}\Bigr(\df{n}{3}\Bigl)\df{q^n}{1-q^n}
+4\sum_{n=1}^{\i}\Bigr(\df{n}{3}\Bigl)\df{q^{4n}}{1-q^{4n}},
\end{equation}
it is easy to show that
\begin{equation}
c(n)=\Bigr( 3-2\delta_{\a,0} \Bigl)(1+(-1)^{\a})\prod_{i=1}^{r}
(1+v_i) \prod_{j=1}^s\df{1+ (-1)^{w_j}}{2},\label{frmcn}
\end{equation}
where $n$ has the prime factorization
\begin{equation*}
n=2^{\a}3^{\b}\prod_{i=1}^{r}p_i^{v_i}\prod_{j=1}^sq_j^{w_j},
\end{equation*}
where  $p_i \equiv 1 \pmod {3}$ and $2\neq q_i \equiv 2 \pmod
{3}$.

From \eqref{tdsss}, we have
\begin{equation}\label{sevdi}
a(n)=\df{c(n)-c(n/3)}{3}+c(n/9)+\df{4}{3}d(n),
\end{equation}
where we assume $c(n/l)=0$ if $l\not|\, n$. If $n \not \equiv 1
\pmod 6$, then $d(n)\equiv 0$ and $a(n)=c(n/9)$ if $3|n$, while
$a(n)=c(n)/3$ if $3\not |n$. This proves the claim in
\eqref{1not6case} for $a(n)$.

Now assume $n \equiv 1 \pmod 6$. If $p$ is a prime and $p \equiv 1
\pmod 3$, then by \eqref{frmcn} we see that $c(p)=4$. If $p$ is
represented by the form $(1,0,27)$, then $a(p)=4$. This is because $4 \leq a(p)
\leq c(p)=4$. Using \eqref{sevdi} we see that $d(p)=2$. If $p$ is not
represented by $(1,0,27)$ then, by \eqref{sevdi}, $d(p)=-1$. Gauss
proved that if $p$ is a prime and $p \equiv 1 \pmod 3$, then $p$ is
represented by $(1,0,27)$ iff $2$ is a cubic residue
modulo $p$ or, equivalently, iff  $2^{(p-1)/3} \equiv 1 \pmod p$.
Therefore,

\begin{equation}\label{dtable}
d(p)=\left\{ \begin{array}{ll}
        2 &  \text{if}\;\; p \equiv 1 \pmod 3 \;\;\text{and}\;\; 2^{(p-1)/3} \equiv 1 \pmod p \\
         -1 &  \text{if}\;\; p \equiv 1 \pmod 3 \;\;\text{and}\;\; 2^{(p-1)/3} \not \equiv 1 \pmod p \\
          0 &\text{if} \;\; p \not \equiv 1 \pmod 3 .\end{array} \right.
\end{equation}

In \cite{martin}, Y.~Martin proved that $qE(q^6)E(q^{18})$ is a
multiplicative cusp form in
$S_1\Bigr(\Gamma_0(108),\Bigr(\df{-108}{*}\Bigl)\Bigl)$. That is
$d(n)$ is multiplicative and for any prime $p$ and $s \geq 0$
\begin{equation}
d(p^{s+2})=d(p)d(p^{s+1})-\Bigr(\df{-108}{p}\Bigl)d(p^s),\label{recur}
\end{equation}
where $d(1)=1$. Using this recursion formula together with
\eqref{dtable}, we find that
\begin{equation}\label{dtable2}
d(p^\a)=\left\{ \begin{array}{ll}
        \a +1 &  \text{if}\;\; p \equiv 1 \pmod 3 \;\;\text{and}\;\; 2^{(p-1)/3} \equiv 1 \pmod p, \\
         \Bigr( (\a +1)| 3 \Bigl) &  \text{if}\;\; p \equiv 1 \pmod 3 \;\;\text{and}\;\; 2^{(p-1)/3} \not \equiv 1 \pmod p, \\
(1+(-1)^{\a})/2 & \text{if} \;\; 2 \neq p \equiv 2 \pmod 3,\\
          0 &\text{if} \;\; p=2 \;\text{or}\;\; 3.\end{array} \right.
\end{equation}
Therefore,
\begin{equation}
d(n)=\delta_{\a,0}\delta_{\b,0}\prod_{i=1}^{r}
(1+v_i)\prod_{i=1}^{s} \Bigr((1+u_i)\mid 3)\Bigl)
\prod_{i=1}^{t}\df{1+ (-1)^{w_i}}{2},\label{frmdn}
\end{equation}
where $n$ has the prime factorization
\begin{equation}
2^{\a}3^{\b}\prod_{i=1}^{r}p_i^{v_i}\prod_{i=1}^{s}q_i^{u_i}\prod_{i=1}^{t}Q_j^{w_i},
\end{equation}
where  $p_i \equiv 1 \pmod {3}$, $2^{\tfrac{p_i-1}{3}} \equiv 1
\pmod {p_i}$,  $ q_i \equiv 1 \pmod {3}$, $2^{\tfrac{q_i-1}{3}}
\not \equiv 1 \pmod {q_i}$ and
 $ 2 \neq Q_i \equiv 2 \pmod {3}$.
By \eqref{sevdi},  if $n \equiv 1 \pmod 6 $, then
$a(n)=\df{c(n)+4d(n)}{3}$. Using \eqref{frmcn} and
\eqref{frmdn} we arrive at the statement for $a(n)$ given in
\eqref{frman}.

From \eqref{lazim1} and \eqref{lazim2}, we have that
\begin{align*}
&\t(q)\t(q^{27})-\sum_{n,m=-\i}^{\i}q^{4n^2+2nm+7m^2}\\&=(\t(q)\t(q^{27})-\t(-q)\t(-q^{27}))/2-2q^7\psi(q^2)\psi(q^{54})\\&=2qE(q^6)E(q^{18}).
\end{align*}

Therefore, $b(n)=a(n)-2d(n)$. The formulas for  $b(n)$ in
\eqref{1not6case} and \eqref{frmbn} now follow from those for
$a(n)$ and $d(n)$. Observe also  from $b(n)=a(n)-2d(n)$ that if
$p$ is a prime and $p \equiv 1 \pmod 3$ then $b(p)=0$ if $p$ is
represented by $(1,0,27)$ and $b(p)=2$, otherwise. Hence, these
primes cannot be represented by $(1,0,27)$ and $(4,2,7)$ at the same
time.
\end{proof}
While they are not explicitly stated there, the formulas for
$a(n)$ and $b(n)$ given by \eqref{1not6case}--\eqref{frmbn} can be
deduced from  Theorem 4.1, Theorem 10.1 and Theorem 10.2 of
\cite{ws} and and Gauss' cubic reciprocity law.

\bigskip

\section{Representations by the Forms $n^2+5m^2+5k^2+5l^2$ and $5n^2+m^2+k^2+l^2$ }

In this section, we give formulas for the number of
representations of positive integers by the quaternary forms
$n^2+5m^2+5k^2+5l^2$ and $5n^2+m^2+k^2+l^2$ and also by the
restricted forms $n+5m+5k+5l$ and $5n+m+k+l$ with $n, m, k,$ and
$l$ being triangular numbers.

\begin{theorem}\label{sno}
\begin{align}
\t(-q)\t^3(-q^5) &=
\Bigr(\df{E^5(q)}{E(q^5)}+4\df{E^5(q^2)}{E(q^{10})}\Bigl)/5 -
\Bigr(q\df{E^5(q^5)}{E(q)}-4q^2\df{E^5(q^{10})}{E(q^2)}
\Bigl),\label{Bidik1} \\
\t^3(-q)\t(-q^5) &=
\Bigr(\df{E^5(q)}{E(q^5)}+4\df{E^5(q^2)}{E(q^{10})}\Bigl)/5 -5
\Bigr(q\df{E^5(q^5)}{E(q)}-4q^2\df{E^5(q^{10})}{E(q^2)}
\Bigl),\label{Bidik2}\\
4q\psi^3(q)\psi(q^5) &=
\Bigr(\df{E^5(q)}{E(q^5)}-\df{E^5(q^2)}{E(q^{10})}\Bigl)/5 +5
\Bigr(q\df{E^5(q^5)}{E(q)}+q^2\df{E^5(q^{10})}{E(q^2)}
\Bigl),\label{Bidik3} \\
4q^2\psi(q)\psi^3(q^5) &=
\Bigr(\df{E^5(q)}{E(q^5)}-\df{E^5(q^2)}{E(q^{10})}\Bigl)/5 +
\Bigr(q\df{E^5(q^5)}{E(q)}+q^2\df{E^5(q^{10})}{E(q^2)}
\Bigl)\label{Bidik4}.
\end{align}
Furthermore, if  $a(n)$, $b(n)$, $c(n)$ and $d(n)$ are defined by
\begin{align*}
\t(q)\t^3(q^5)=:1+\sum_{n=1}^{\i}a(n)q^n, \;\;\;\;
\t^3(q)\t(q^5)=:1+\sum_{n=1}^{\i}b(n)q^n,\\
4q\psi^3(q)\psi(q^5)=:\sum_{n=1}^{\i}c(n)q^n,\;\;\;\;
4q^2\psi(q)\psi^3(q^5)=:\sum_{n=1}^{\i}d(n)q^n,
\end{align*}
then
\begin{align}
a(n)&= (-1)^{n-1}(1+5^d(-1)^{g+t}) \df{
(5+(-2)^{g+1})}{3}\prod_{i=1}^{r} \df{1-p_i^{v_i+1}}{ 1-p_i}
\prod_{j=1}^{s}\df{1- (-q_j)^{w_j+1}}{1+q_j},\label{form1}\\
b(n)&= (-1)^{n-1}(1+5^{d+1}(-1)^{g+t}) \df{
(5+(-2)^{g+1})}{3}\prod_{i=1}^{r} \df{1-p_i^{v_i+1}}{ 1-p_i}
\prod_{j=1}^{s}\df{1- (-q_j)^{w_j+1}}{1+q_j},\label{form2}\\
c(n)&= (-2)^g(-1+5^{d+1}(-1)^{g+t}) \prod_{i=1}^{r}
\df{1-p_i^{v_i+1}}{ 1-p_i}
\prod_{j=1}^{s}\df{1- (-q_j)^{w_j+1}}{1+q_j},\label{form3}\\
d(n)&= (-2)^g(-1+5^d(-1)^{g+t})\prod_{i=1}^{r} \df{1-p_i^{v_i+1}}{
1-p_i} \prod_{j=1}^{s}\df{1- (-q_j)^{w_j+1}}{1+q_j},\label{form4}
\end{align}
where $n$ has the prime factorization
\begin{equation*}
n=2^g5^d\prod_{i=1}^{r}p_i^{v_i}\prod_{j=1}^sq_j^{w_j},
\end{equation*}
with  $p_i \equiv \pm 1\pmod {5}$ and $q_i \equiv \pm 2 \pmod
{5}$, $q_i$ is odd and $t$ is the number of odd prime divisors of
$n$, counting multiplicities, that are congruent to $\pm  2
\pmod{5}$.
\end{theorem}

\begin{proof}
Our proof employs the following well known Lambert series
identities of Ramanujan

\begin{align}
\df{E^5(q)}{E(q^5)}&=1-5\sum_{n=1}^{\i}\Bigr(\df{n}{5}\Bigl)\df{nq^n}{1-q^n},\label{c5core}\\
q\df{E^5(q^5)}{E(q)}&=\sum_{n=1}^{\i}\Bigr(\df{n}{5}\Bigl)\df{q^n}{(1-q^n)^2}.\label{5core}
\end{align}
For the history of these and many related identities see \cite{GAn},
\cite[pp.~249--263]{III}. We also use  the theta function
identities \cite[p.~262, Entry 10]{III}
\begin{equation}\label{ttid}
\t^2(q)-\t^2(q^5)=4qf(q,q^9)f(q^3,q^7)
\end{equation}
and
\begin{equation}\label{ssid}
\psi^2(q)-q\psi^2(q^5)=f(q^2,q^3)f(q,q^4).
\end{equation}

By multiplying both sides of \eqref{ttid} with $\t^3(q^5)/ \t(q)$,
we find that
\begin{equation}
\t(q)\t^3(q^5)-\df{\t^5(q^5)}{\t(q)}=4q\df{E^5(-q^5)}{E(-q)}.\label{bidik1}
\end{equation}
From \eqref{bidik1}, we deduce that
\begin{align}
16q^2\df{E^5(q^{10})}{E(q^2)}&=16q^2\df{\t(q)}{\t^5(q^5)}\df{E^{10}(-q^5)}{E^2(-q)}\notag\\
&=\df{\t(q)}{\t^5(q^5)}\Bigr\{\t(q)\t^3(q^5)-\df{\t^5(q^5)}{\t(q)}
\Bigl\}^2\notag\\
&=\t^3(q)\t(q^5)-2\t(q)\t^3(q^5)+\df{\t^5(q^5)}{\t(q)}.\label{bidik2}
\end{align}
Using the imaginary transformation on \eqref{bidik1} and
\eqref{bidik2}, we obtain, respectively, that
\begin{equation}
5\t^3(q)\t(q^5)-\df{\t^5(q)}{\t(q^5)}=4\df{E^5(-q)}{E(-q^5)}\label{bidik3}
\end{equation}
and
\begin{equation}
25\t(q)\t^3(q^5)-10\t^3(q)\t(q^5)+\df{\t^5(q)}{\t(q^5)}=16\df{E^5(q^2)}{E(q^{10})}.\label{bidik4}
\end{equation}

By multiplying both sides of \eqref{ssid} with $\psi^3(q^5)/
\psi(q)$, we find that
\begin{equation}
\psi(q)\psi^3(q^5)-q\df{\psi^5(q^5)}{\psi(q)}=\df{E^5(q^{10})}{E(q^2)}.\label{2bidik1}
\end{equation}
From \eqref{2bidik1}, we also find that
\begin{align}
\df{E^5(q^{5})}{E(q)}&=\df{\psi(q)}{\psi^5(q^5)}\df{E^{10}(q^{10})}{E^2(q^2)}\notag\\
&=\df{\psi(q)}{\psi^5(q^5)}\Bigr\{\psi(q)\psi^3(q^5)-q\df{\psi^5(q^5)}{\psi(q)}
\Bigl\}^2\notag\\
&=\psi^3(q)\psi(q^5)-2q\psi(q)\psi^3(q^5)+q^2\df{\psi^5(q^5)}{\psi(q)}.\label{2bidik2}
\end{align}

Using the imaginary transformation on \eqref{2bidik1} and
\eqref{2bidik2}, we obtain, respectively, that
\begin{equation}
-5q\psi^3(q)\psi(q^5)+\df{\psi^5(q)}{\psi(q^5)}=\df{E^5(q^{2})}{E(q^{10})}\label{2bidik3}
\end{equation}
and
\begin{equation}
25q^2\psi(q)\psi^3(q^5)-10q\psi^3(q)\psi(q^5)+\df{\psi^5(q)}{\psi(q^5)}=\df{E^5(q)}{E(q^5)}.\label{2bidik4}
\end{equation}
Using \eqref{bidik1}--\eqref{2bidik4}, we easily derive
\eqref{Bidik1}--\eqref{Bidik4}.

Next, we sketch a proof of \eqref{form1}. We omit the proofs of
\eqref{form2}--\eqref{form4} since their proofs are similar to
that of \eqref{form1}. For convenience, $[q^n]V(q)$ will denote
the coefficient of $q^n$ in the Taylor series expansion of $V(q)$.\\
From \eqref{c5core}, we have

\begin{equation}
\df{E^5(q)}{E(q^5)}=1-5\sum_{n=1}^{\i}\Bigr(\df{n}{5}\Bigl)\df{nq^n}{1-q^n}=1-5\sum_{n=1}^{\i}\Bigr(\sum_{d|n}\Bigr(\df{d}{5}\Bigl)d\Bigl)q^n.
\end{equation}
Using the fact that the coefficients are given by the
multiplicative function $\sum_{d|n}\Bigr(\df{d}{5}\Bigl)d$, we
conclude that for $n>1$
\begin{equation}
[q^n]\df{E^5(q)}{E(q^5)}=-5\prod_{i=1}^{r} \df{1-p_i^{v_i+1}}{
1-p_i} \prod_{j=1}^s\df{1-  (-q_j)^{w_j+1}}{1+q_j},\label{xz1}
\end{equation}
where $n$ has the prime factorization
\begin{equation}
n=5^d\prod_{i=1}^{r}p_i^{v_i}\prod_{j=1}^sq_j^{w_j},
\end{equation}
with  $p_i \equiv \pm 1\pmod {5}$ and $q_i \equiv \pm 2 \pmod
{5}$.

It is easy to show \cite[Thr.~4]{GKS}
\begin{equation}
[q^n]\df{qE^5(q^5)}{E(q)}=5^d\prod_{i=1}^{r} \df{1-p_i^{v_i+1}}{
1-p_i} \prod_{j=1}^{s}(-1)^{w_j}\df{1-
(-q_j)^{w_j+1}}{1+q_j},\label{xz2}
\end{equation}
where $n>0$ has the prime factorization
\begin{equation*}
n=5^d\prod_{i=1}^{r}p_i^{v_i}\prod_{j=1}^sq_j^{w_j},
\end{equation*}
with  $p_i \equiv \pm 1\pmod {5}$ and $q_i \equiv \pm 2 \pmod
{5}$. Using \eqref{xz1}, \eqref{xz2} together with \eqref{Bidik1},
we arrive at \eqref{form1}.
\end{proof}

Theorem \ref{sno} has the following
\begin{corollary}
\begin{align*}
&a) \;\;\;\; [q^n] (\t^3(q)\t(q^5)) >0 , \;\;[q^n]
(\psi^3(q)\psi(q^5))
>0
\;\;\text{for any}\; n \geq 0,\\
&b) \;\;\;\;[q^n] (\t(q)\t^3(q^5)) =0 , \;\;[q^n]
(\psi(q)\psi^3(q^5))
>0 \;\;\text{iff}\;\; n \equiv 2\; \text{or}\; 3 \pmod 5.
\end{align*}
\end{corollary}
Note that our  corollary is in  agreement   with Ramanujan's
observation in  \cite{rmn}  where  quadratic  form $x^2+y^2+z^2+ 5
w^2$ is listed  as universal. It means that this form represents all positive
integers. Interested reader may want to check  \cite{bgv} for
new results  about  universal quadratic forms.

Next,  we use \eqref{bidik1}--\eqref{bidik3} to derive
\begin{align}
&\Bigr(\df{E^5(-q)}{E(-q^5)}+4\df{E^5(q^2)}{E(q^{10})}\Bigl)/5\notag\\
&=(5\t(q)\t^3(q^5)-\t^3(q)\t(q^5))/4\\
&=\df{1}{4}\df{\t^2(q^5)}{\t^2(q)}\Bigr\{
5\t^3(q)\t(q^5)-\df{\t^5(q)}{\t(q^5)}\Bigl\}\\
&=\df{\t^2(q^5)E^5(-q)}{\t^2(q)E(-q^5)}\\
&=\df{E^5(q^2)E^7(q^{10})}{E(q)E(q^4)E^3(q^5)E^3(q^{20})}.\label{88p1}
\end{align}
Moreover, using \eqref{c5core} and \eqref{5core} we obtain
\begin{align}
&\Bigr(\df{E^5(-q)}{E(-q^5)}+4\df{E^5(q^2)}{E(q^{10})}\Bigl)/5\notag\\
&=1+\sum_{n=1}^{\i}\Bigr(\df{n}{5}\Bigl)\df{nq^n}{1-(-q)^n}\label{nnls1}\\
&=1+\sum_{n=1}^{\i}\Bigr(\sum_{d|n}(-1)^{n+d}d\Bigr(\df{d}{5}\Bigl)\Bigl)q^n.
\end{align}
Arguing as above one finds
\begin{align}
q\df{E^5(-q^5)}{E(-q)}+4q^2\df{E^5(q^{10})}{E(q^2)}&=q\df{\t^2(q)E^5(-q^5)}{\t^2(q^5)E(-q)}\label{sdd19}\\
&=q\df{E^7(q^2)E^5(q^{10})}{E^3(q)E^3(q^4)E(q^5)E(q^{20})}\label{88p2}\\
&=-\sum_{n=1}^{\i}\Bigr(\df{n}{5}\Bigl)\df{(-q)^n}{(1+(-q)^n)^2}\label{nnnls2}\\
&=\sum_{n=1}^{\i}\Bigr(\sum_{d|n}(-1)^{n+d}d\Bigr(\df{n/d}{5}\Bigl)\Bigl)q^n,
\end{align}
and
\begin{align}
-\Bigr(\df{E^5(q)}{E(q^5)}-\df{E^5(q^2)}{E(q^{10})}\Bigl)/5 &=q\df{\psi^2(q^5)E^5(q^2)}{\psi^2(q)E(q^{10})}\\
&=q\df{E^2(q)E(q^2)E^3(q^{10})}{E^2(q^5)}\label{88p3}\\
&=\sum_{n=1}^{\i}\Bigr(\df{n}{5}\Bigl)\df{nq^n}{1-q^{2n}}\label{nnls2}\\
&=\sum_{n=1}^{\i}\Bigr(\sum_{d|n}d\Bigr(\df{d}{5}\Bigl)\gamma(n/d)\Bigl)q^n\label{Lazimm3},
\end{align}
where
\begin{equation*}
\gamma(n):=\left\{ \begin{array}{ll}
        1 &  \text{if $n$ is odd}, \\
        0 &  \text{if $n$ is even} .\end{array} \right.
\end{equation*}
Lastly,
\begin{align}
q\df{E^5(q^5)}{E(q)}+q^2\df{E^5(q^{10})}{E(q^2)}
&=q\df{\psi^2(q)E^5(q^{10})}{\psi^2(q^{5})E(q^2)}\label{sdd20}\\
&=q\df{E^3(q^2)E^2(q^5)E(q^{10})}{E^2(q)}\label{88p4}\\
&=\sum_{\substack{n=1\\n\, \text{is odd}}}^{\i}\Bigr(\df{n}{5}\Bigl)\df{q^n}{(1-q^{n})^2}\label{nnnls3}\\
&=\sum_{n=1}^{\i}\Bigr(\sum_{d|n}\gamma(d)\Bigr(\df{d}{5}\Bigl)n/d\Bigl)q^n.
\end{align}
The $eta$-quotients given by \eqref{88p3} and \eqref{88p4} are
clearly multiplicative. Using the fact that $(-1)^{n+1}$ is a
multiplicative function of $n$, we also see that the first two
$eta$-quotients given in \eqref{88p1} and \eqref{88p2} are
multiplicative. These four multiplicative $eta$-quotients are not
included in Martin's list. Lambert series representations of
\eqref{nnls1} and \eqref{nnls2} are given by Ramanujan
\cite[p.~249, Entry 8 (i), (ii)]{III}. The identity
\eqref{sdd20}--\eqref{88p4} is the identity (6.12) of \cite{BG}.
Using \eqref{Bidik1}--\eqref{Bidik4}, \eqref{nnls1},
\eqref{nnnls2}, \eqref{nnls2} and \eqref{nnnls3}, we have
\begin{corollary}
\begin{align}
\t(q)\t^3(q^5)&=1+\sum_{n=1}^{\i}\Bigr(\df{n}{5}\Bigl)
\df{nq^n}{1-(-q)^n}-\sum_{n=1}^{\i}\Bigr(\df{n}{5}\Bigl)\df{(-q)^n}{(1+(-q)^n)^2},\\
\t^3(q)\t(q^5)&=1+\sum_{n=1}^{\i}\Bigr(\df{n}{5}\Bigl)
\df{nq^n}{1-(-q)^n}-5\sum_{n=1}^{\i}\Bigr(\df{n}{5}\Bigl)\df{(-q)^n}{(1+(-q)^n)^2},\\
4q\psi^3(q)\psi(q^5)&=-\sum_{n=1}^{\i}\Bigr(\df{n}{5}\Bigl)\df{nq^n}{1-q^{2n}}+5\sum_{\substack{n=1\\n\,\text{is odd}}}^{\i}\Bigr(\df{n}{5}\Bigl)\df{q^n}{(1-q^{n})^2},\\
4q^2\psi(q)\psi^3(q^5)&=-\sum_{n=1}^{\i}\Bigr(\df{n}{5}\Bigl)\df{nq^n}{1-q^{2n}}+\sum_{\substack{n=1\\n\,\text{is
odd}}}^{\i}\Bigr(\df{n}{5}\Bigl)\df{q^n}{(1-q^{n})^2}.
\end{align}
\end{corollary}

\section{Outlook}
Clearly,  this manuscript does not  exhaust all potential
connections  between the Ramanujan identities  and quadratic
forms. We believe that  Ramanujan identities can be employed  to
find coefficients of  many sextenary forms. For example,  in
\cite{BHBH} we will show how to use  our new identity
\begin{align*}
7\t^3(-q)\t^3(-q^7)
=& -49\Bigr(q^2\df{E^7(q^7)}{E(q)}+qE^3(q)E^3(q^7)\Bigl)\\
&+56\Bigr(7q^4\df{E^7(q^{14})}{E(q^2)}+q^2E^3(q^2)E^3(q^{14})\Bigl)
-\df{E^7(q)}{E(q^7)}+8\df{E^7(q^2)}{E(q^{14})},
\end{align*}
together with two identities of Ramanujan to determine the
coefficients of $\t^3(q)\t^3(q^7)$. There   we   shall also prove
the following intriguing inequalities
\begin{align*}
[q^n] \Bigr(\psi^3(q)\psi^3(q^7)  -q\df{E^7(q^{14})}{E(q^2)}\Bigl) \geq 0, \\
[q^n] \Bigr( \t^3(q)\t^3(q^7) + q^7  -q^2\df{E^7(q^7)}{E(q)}\Bigl)
\geq 0.
\end{align*}
We will also determine the coefficients of $\t^5(q)\t(q^3)$ and
$\t(q)\t^5(q^3)$ by proving that
\begin{align*}
\t^5(q)\t(q^3)=1+\sum_{n=1}^{\i}\Bigr(\sum_{d|n}(-1)^{n+d}d^2\Bigr(\df{d}{3}\Bigl)\Bigl)q^n
+9\sum_{n=1}^{\i}\Bigr(\sum_{d|n}(-1)^{n+d}d^2\Bigr(\df{n/d}{3}\Bigl)\Bigl)q^n
\end{align*}
and
\begin{align*}
\t(q)\t^5(q^3)=1+\sum_{n=1}^{\i}\Bigr(\sum_{d|n}(-1)^{n+d}d^2\Bigr(\df{d}{3}\Bigl)\Bigl)q^n
+\sum_{n=1}^{\i}\Bigr(\sum_{d|n}(-1)^{n+d}d^2\Bigr(\df{n/d}{3}\Bigl)\Bigl)q^n.
\end{align*}

In the course of this investigation we have determined the
coefficients of many multiplicative  $eta$-quotients. Those
$eta$-quotients given by the equations \eqref{d20p1},
\eqref{d20p2}, \eqref{p1p2def}, and \eqref{lazim2} are on Martin's
list \cite{martin}, while those that are given by \eqref{fix60i},
\eqref{fix60ii}, \eqref{fix60iii}, \eqref{88p1}, \eqref{88p2},
\eqref{88p3} and \eqref{88p4} are not in his list. This is because
Martin only considered  multiplicative  $eta$-quotients which are
eigenforms for all Hecke operators. As an example, it seems that
the multiplicative $eta$-quotient  $-Q(-q)$,  defined by
\eqref{fix60ii}, is an eigenform for all Hecke operators $T_p$
with odd prime $p$ but it is not hard to show that
\begin{equation*}
T_2(-Q(-q))=\df{E^2(q^6)E^2(q^{10})}{E(q^2)E(q^{30})}.
\end{equation*}

The $eta$-quotients in \eqref{c5core} and \eqref{5core} also
appear in Martin's list \cite{martin}. In our future publications
we plan to discuss the coefficients of all multiplicative
$eta$-products that appear on this list. To this end we proved the
following

\begin{theorem}
Suppose the prime factorization of $n$ is
\begin{equation*}
2^{\a}\prod_{i=1}^{r}p_i^{u_i}\prod_{i=1}^sq_i^{v_j}\prod_{i=1}^{l}P_i^{w_i}\prod_{i=1}^{t}Q_j^{d_i},
\end{equation*}
where
\begin{align*}
Q_i &\equiv 3 \pmod 4,\\
p_i &\equiv 5 \pmod 8, \\
q_i &\equiv 1 \pmod 8 \;\; \text{and}\;\; 2^ {(q_i-1)/4}  \equiv
1 \pmod {q_i},\\
P_i &\equiv 1 \pmod 8 \;\; \text{and}\;\;  2^ {(P_i-1)/4}  \equiv
-1 \pmod {P_i}.
\end{align*}

Then
\begin{equation*}
[q^n]\Bigr( q\df{ E^4(q^{16})}{ E(q^{32}) E(q^8)}\Bigl)
=\delta_{\a,0}\prod_{i=1}^{r}(-1)^{u_i}(1+(-1)^{u_i})/2\prod_{i=1}^{s}(1+v_i)\prod_{i=1}^l(-1)^{w_j}(1+w_j)\prod_{i=1}^{t}(1+(-1)^{d_i})/2.
\end{equation*}
\end{theorem}
We  would like to conclude with a following remarkable identity:
\begin{equation*}
q\df{\widetilde{Q}(2,0,7)+
\widetilde{Q}(3,2,5)}{\widetilde{Q}(1,0,14)- \widetilde{Q}(3,2,5)}
=\df{E^2(q^{4})E^2(q^{14})}{ E^2(q^{2})E^2(q^{28})},
\end{equation*}
where
\begin{equation*}
\widetilde{Q}(a,b,c):=\sum_{n,m=-\i}^{\i} q^{an^2+ bnm+cm^2}.
\end{equation*}
It is important to observe that  neither $\widetilde{Q}(2,0,7)+
\widetilde{Q}(3,2,5)$ nor $\widetilde{Q}(1,0,14)-
\widetilde{Q}(3,2,5)$ is an $eta$-quotient. This result along with
other similar type identities will be discussed elsewhere.

\end{document}